\documentclass[journal]{IEEEtran}
\usepackage{etex}
\reserveinserts{100}
\usepackage{mathtools}
\usepackage{graphics, theorem, times, amsfonts, graphicx, amssymb, cite}
\usepackage{tikz}
\usetikzlibrary{shapes,arrows}
\usepackage{pgfplots}
\usepackage{color}
\usepackage{setspace}
\usepackage{hyperref}
\usepackage{multirow}
\usepackage{rotating}
\usepackage{comment}
\usepackage{flushend}
\usepackage[font=small]{caption}
\usepackage{subcaption}
\usepackage[affil-it]{authblk}
\usepackage{bm}
\usepackage{algorithm}
\usepackage{algorithmic}
\usepackage{enumerate}
\usepackage{morefloats}
\input{mysymbol.sty}
\usepackage{needspace}

\renewcommand{\baselinestretch}{.995}
\newcommand{\QED}{\hfill\ensuremath{\blacksquare}}

\newtheorem{assumption}{\hspace{0pt}\bf Assumption}
\newtheorem{lemma}{\hspace{0pt}\bf Lemma}
\newtheorem{proposition}{\hspace{0pt}\bf Proposition}

\newtheorem{theorem}{\hspace{0pt}\bf Theorem}

\newtheorem{remark}{\hspace{0pt}\bf Remark}

\title{Decentralized Quasi-Newton Methods}
\author{Mark Eisen, Aryan Mokhtari, and Alejandro Ribeiro 
\thanks{{Supported by NSF CAREER CCF-0952867 and ONR N00014-12-1-0997. The authors are with the Department of Electrical and Systems Engineering, University of Pennsylvania. Email at: \{maeisen, aryanm, aribeiro\}@seas.upenn.edu.
}}}

\begin{document}

\thispagestyle{empty}
\maketitle

\begin{abstract}
We introduce the decentralized Broyden-Fletcher-Goldfarb-Shanno (D-BFGS) method as a variation of the BFGS quasi-Newton method for solving decentralized optimization problems. The D-BFGS method is of interest in problems that are not well conditioned, making first order decentralized methods ineffective, and in which second order information is not readily available, making second order decentralized methods impossible. D-BFGS is a fully distributed algorithm in which nodes approximate curvature information of themselves and their neighbors through the satisfaction of a secant condition. We additionally provide a formulation of the algorithm in asynchronous settings. Convergence of D-BFGS is established formally in both the synchronous and asynchronous settings and strong performance advantages relative to first order methods are shown numerically.
\end{abstract}


\section{Introduction} \label{sec_intro}

Decentralized optimization involves a group of interconnected agents seeking to jointly minimize a common objective function about which they have access to information that is local and partial. The agents collaborate by successively sharing information with other agents located in their communication neighborhood with the goal of eventually converging to the network-wide optimal argument. Decentralized optimization has proven effective in contexts where information is gathered by different nodes of a network, such as decentralized control
\cite{Bullo2009,Cao2013-TII,LopesEtal8}, wireless systems
\cite{Ribeiro10,scutari2014decomposition,Ribeiro12}, sensor networks
\cite{Schizas2008-1,KhanEtal10,cRabbatNowak04,barbarossa2007decentralized}, and large scale
machine learning
\cite{bekkerman2011scaling,Tsianos2012-allerton-consensus,Cevher2014}.

Although there are different formulations of decentralized optimization problems, all have in common a reliance on the distributed computability of the gradient. This property refers to the ability of each agent to compute gradients with respect to its local variable using its own variable and the variables of neighboring nodes. If this property holds, it is possible for nodes to exchange variables with neighbors, compute gradients with respect to their local variables, implement the corresponding block of a gradient descent algorithm, and proceed to a new variable exchange to repeat the process. Distributed gradient computability is sometimes inherent to the objective function \cite{eksin2012distributed}, but more often the result of some reformulation. The latter is the case in consensus optimization problems which do not have distributedly computable gradients but can be transformed into problems where the gradients are (see Section \ref{sec_consensus}). The most popular techniques for doing so are the use of penalties to enforce the consensus constraint \cite{nedic2009, YuanQing, Jakovetic2014-1,NN-part1} and the use of gradient ascent in the dual domain for a problem in which consensus is imposed as a constraint {\cite{cRabbatNowak04,Schizas2008-1,makhdoumi2016convergence,bianchi2014stochastic}}. 

The problem with methods that rely on distributed gradient computations is that gradient descent methods exhibit slow convergence. This limits applicability to cases where the function to be optimized is well conditioned, which in practice implies arguments with low dimension. The Hessian can be used to determine a better descent direction if it so happens to also be distributedly computable -- which it is if the Hessian matrix has the same sparsity pattern of the network. This is not to say that Newton's method can be implemented in a distributed manner, because the non-sparse Hessian inverse would be needed for that. Still, the Hessian can be used to approximate the Newton step and yield methods with faster convergence. This has been done for consensus optimization problems reformulated as penalty methods \cite{NN-part1,NN-part2} and for the dual problem of optimal linear flow control \cite{zargham2014accelerated}. These approximate Newton methods exhibit faster convergence relative to their corresponding first order methods.

An alternative to the approximation of Newton steps is the use of quasi-Newton methods that rely on gradients to produce a curvature estimation to use in lieu of the Hessian inverse \cite{dennis1974characterization,powell1976some}. The goal of this paper is to adapt the curvature estimation technique of the Broyden-Fletcher-Goldfarb-Shanno (BFGS) quasi-Newton optimization method for use in distributed settings. This adaptation leads to the development of the distributed (D)-BFGS method that we show can be implemented with nodes that operate either synchronously or asynchronously. We further prove convergence for both methods in the case of convex function and establish a linear convergence rate for the case of strongly convex functions in both, the synchronous and the asynchronous formulations. The advantages of D-BFGS relative to approximate Newton methods are that they do not require computation of Hessians, which can itself be expensive, and that they apply in any scenario in which gradients are distributedly computable irrespectively of the structure of the Hessian. 

The paper starts with the introduction of notation and a formal statement of the gradient distributed computability property (Section \ref{sec_problem_formulation}). The curvature approximation used in BFGS is then introduced (Section \ref{sec_grad_descent}). The fundamental observation here is that this curvature approximation is chosen to satisfy a secant condition because this is a property that the Hessian has. We then point out that the secant condition can be written distributedly as long as gradients are distributedly computable. Building on this observation we define D-BFGS as a method where the Hessian inverse is approximated by a matrix that satisfies the secant condition but whose sparsity pattern is chosen a fortiori to match the sparsity pattern of the graph (Section \ref{sec_dbfgs}). This matching of sparsity patterns guarantees that the method can be implemented in a distributed manner (Algorithm \ref{alg_dbfgs}). The D-BFGS method requires three separate variable exchanges in each iteration. Since the time cost of this synchronization can be significant, we introduce an asynchronous version where nodes operate on their local memories which are synchronized by a communication protocol that runs on a separate clock (Section \ref{sec_async_dbfgs}). In this asynchronous algorithm nodes operate with possibly -- indeed, most often -- outdated variables to avoid the time cost of running in synch (Algorithm \ref{alg_async_dbfgs}).

Convergence properties are then established (Section \ref{sec_convergence}). In the case of synchronous D-BFGS we prove convergence for smooth convex functions and further establish a linear rate when the functions are strongly convex (Section \ref{sec_synch_convergence}). For the case of asynchronous implementations we impose an upper bound in the number of iterations that it takes for the information of a node to get updated in the local memory of its neighbors. Under this hypothesis we also establish convergence for smooth convex functions and a linear rate for strongly convex functions (Section \ref{sec_asynch_convergence}). The convergence rate decreases with increasing levels of asynchronicity. The application of D-BFGS in consensus optimization problems is then explicitly discussed (Section \ref{sec_consensus}). We explain how D-BFGS can be used in combination with penalty methods to obtain a quasi-Newton version of distributed gradient descent (Section \ref{sec_primal_consensus}) and how it can be used in the dual domain to obtain a quasi-Newton version of distributed dual ascent (Section \ref{sec_dual_consensus}). We close the paper with numerical results comparing the performance of D-BFGS to first order methods on various consensus problems in both the synchronous and asynchronous settings (Section \ref{sec_numerical_results}).

\section{Problem Formulation} \label{sec_problem_formulation}

Consider a decentralized system of $n$ nodes, each of which has access to a local variable $\bbx_i \in \reals^{p}$. Nodes are connected by a communication graph $\ccalG=(\ccalV,\ccalE)$ with nodes $\ccalV=\{1,\dots,n\}$ and $m$ edges $\ccalE=\{(i,j)\ |\ i\ \text{and}\ j \ \text{are connected} \}$. We assume the graph $\ccalG$ is undirected which implies $(i,j)\in\ccalE$ if and only if $(j,i)\in \ccalE $. Define the set $n_i$ as the neighborhood of node $i$ including $i$, i.e., $n_i=\{j\ |\ j=i \lor (i,j)\in\ccalE\}$, and the neighborhood size $m_i := |n_i|$. Further define the global variable $\bbx=[\bbx_1;\dots;\bbx_n]\in\reals^{np}$ as the concatenation of the local variables $\bbx_i$ and, for each $i$, the neighborhood variable $\bbx_{n_i}=\{ \bbx_j \}_{j \in n_i} \in\reals^{m_i p}$ as the concatenation of local variables belonging to nodes in $n_i$. The system's goal is to find the optimal argument $\bbx^*  \in \reals^{np}$ that minimizes the smooth convex cost function  $f: \reals^{np} \rightarrow \reals$ when the gradient components $\nabla_i f(\bbx):=\partial f(\bbx)/\partial \bbx_i$ have a local structure,
\begin{align}\label{eq_gen_problem}
   \bbx^* \ := & \ \argmin_{\bbx \in \reals^{np}} \ f(\bbx), \quad 
                 \text{ with\ } \nabla_i f(\bbx) = \nabla_i f(\bbx_{n_i}) .
\end{align}
Since they are functions of variables that are available in their respective neighborhoods, gradient components $ \nabla_i f(\bbx)$ can be evaluated at node $i$ using only single hop communications. We study examples of network optimization problems with gradients that have this property in Section \ref{sec_consensus}.

\subsection{Gradient descent and BFGS} \label{sec_grad_descent}

The gradient property in \eqref{eq_gen_problem} means that it is possible to implement gradient descent on $f(\bbx)$  in a distributed manner whereby the $i$th component of $\bbx$ is updated iteratively at node $i$ until it converges to the $i$th component of the optimal solution $\bbx^*$. Introduce then the time index $t$ and the variable $\bbx(t)$ to be its value at time $t$ and define the update
\begin{equation}\label{eq_descent_update}
   \bbx(t+1) = \bbx(t) + \epsilon(t) \bbd(t),
\end{equation}
where $\epsilon(t)$ is a scalar stepsize. For convex functions, convergence of $\bbx(t)$ to $\bbx^*$ is guaranteed if $\bbd(t)$ is a proper descent direction for which $\bbd(t)^T \nabla f(\bbx(t))\leq0$. Since the negative gradient has this property \cite{boyd}, a natural choice is to make
\begin{equation}\label{eq_grad_descent}
   \bbd(t)  = - \nabla f(\bbx(t)) : = -\bbg(t).
\end{equation}
Using the descent direction in \eqref{eq_grad_descent} for the update in \eqref{eq_descent_update} yields the gradient descent method \cite{boyd}. The corresponding iterations can be written componentwise as $\bbx_i(t+1) = \bbx_i(t) + \epsilon(t) \nabla_i f(\bbx(t))$ and can be implemented in a distributed manner if the gradient $\nabla_i f(\bbx(t))$ satisfies the property in \eqref{eq_gen_problem}.

As in the centralized case, decentralized gradient descent methods are often slow to converge. In centralized systems, speed of convergence can be increased by premultiplying $\bbg(t)$ by a positive definite matrix to obtain a better descent direction. Newton's method premultiplies the gradient with the Hessian inverse $\nabla^2 f(\bbx)^{-1}$ and is recommended whenever possible \cite{NN-part1, NN-part2}. However, the use of Hessian inverses in a distributed implementation requires further assumptions on \eqref{eq_gen_problem} specific to the problem of interest and leads to problem specific challenges and limitations; see, e.g., \cite{NN-part1}.  
 
Alternatively, quasi-Newton methods approximate the objective function Hessian inverse using subsequent gradient evaluations. To be more precise, define the descent direction 
\begin{equation}\label{eq_qnewt_update}
\bbd(t) = -\bbB(t)^{-1}\bbg(t),
\end{equation}
where $\bbB(t)$ is a symmetric positive definite matrix that serves as an approximation of the Hessian $\nabla^2 f(\bbx(t))$. Various quasi-Newton methods differ in how they define $\bbB(t)$, with the most common being the method of Broyden-Fletcher-Goldfarb-Shanno (BFGS). To formulate BFGS begin by defining the variable variation $\bbv(t)$ and the gradient variation $\bbr(t)$ vectors, 
\begin{align}
\bbv(t) = \bbx(t+1) - \bbx(t), \qquad \bbr(t) = \bbg(t+1) - \bbg(t).\label{eq_bfgs_vars} 
\end{align}
Observe that $\bbv(t)$ and $\bbr(t)$ capture differences of two consecutive variables and gradients, respectively, evaluated at steps $t+1$ and $t$. At each iteration, we select a new Hessian approximation $\bbB(t+1)$ that satisfies the secant condition $\bbB(t+1) \bbv(t) = \bbr(t)$. This condition is fundamental, as the secant condition is satisfied by the actual Hessian for small $\bbv(t)$. As this is an underdetermined system, we select $\bbB(t+1)$ such that it is closest to the previous approximation in terms of Gaussian differential entropy,
\begin{alignat}{2}\label{eq_bfgs_update}
\bbB(t+1) = &\argmin_{\bbZ}\ &&
\text{tr}[ \bbB(t)^{-1} \bbZ] - \text{logdet}[\bbB(t)^{-1} \bbZ] - n, \nonumber\\
&  \text{   s.t.}\quad &&
\bbZ \bbv(t) = \bbr(t),\quad \bbZ \succeq \bb0.
\end{alignat}
Note that we also require the next approximation to be positive semidefinite to ensure a proper descent. In order for the problem to be feasible with a positive definite solution, it is necessary to have $\bbv(t)^T \bbr(t) > 0$. This is always true when the objective function is strongly convex \cite{mokhtari2014res}. The closed-form solution to \eqref{eq_bfgs_update} provides the BFGS update formula
\begin{align}\label{eq_bfgs_cent} 
\bbB(t+1) &=  \bbB(t) + \frac{\bbr(t) \bbr(t)^T}{\bbr(t)^T \bbv(t)} - \frac{\bbB(t) \bbv(t) \bbv(t)^T \bbB(t)}{\bbv(t)^T \bbB(t) \bbv(t)},
\end{align}
which shows that $\bbB(t+1)$ can be computed using the previous approximation matrix $\bbB(t)$ as well as the variable $\bbv(t)$ and gradient $\bbr(t)$ variations at step $t$.

The matrices $\bbB(t)$ that solve \eqref{eq_bfgs_update} -- which are explicitly given in \eqref{eq_bfgs_cent} -- depend only on gradients that we are assuming can be computed in a distributed manner. This does {\it not} mean that a distributed implementation of BFGS is possible because neither $\bbB(t)$ nor $\bbB(t)^{-1}$ have a sparsity pattern to permit local evaluation of descent directions. Additionally, the computation of $\bbB(t+1)$ in \eqref{eq_bfgs_cent} includes the inner product $\bbr(t)^T \bbv(t)$, which itself requires global information. It is important to note, however, is that the choice of objective function in \eqref{eq_bfgs_update} is of secondary importance to satisfying the secant condition and the secant condition {\it does} have a structure that allows for distributed evaluation. In the following section we resolve the issue of decentralization by introducing a variation of BFGS which modifies the Hessian approximation such that $\bbd(t)$ is computable distributedly.

\medskip\noindent{\bf Notation remark. } The $i$th block of a vector $\bbz \in \reals^{np}$ is denoted as $\bbz_i \in \reals^{p}$, while $\bbz_{n_i} \in \reals^{m_i p}$ denotes the components in $n_i$. To have global representations, we define $\hbz_{n_i} \in \reals^{np}$ to be the vector $\bbz_{n_i}$ padded with zeros in locations corresponding to nodes not in $n_i$. Likewise, for any matrix $\bbA \in \reals^{np \times np}$, we define $\bbA_{n_i} \in \reals^{m_i p \times m_i p}$ to be the $m_i p$ rows and columns of $\bbA$ corresponding to nodes in $n_i$ and $\hbA_{n_i} \in \reals^{n p \times n p}$ to be the matrix $\bbA_{n_i}$ padded with zeros in other locations.

%
\section{Decentralized BFGS}\label{sec_dbfgs}

Our goal here is to develop an algorithm of the form
\begin{equation}\label{eq_descent_update_local}
\bbx_i(t+1) = \bbx_i(t) + \epsilon(t) \bbd_i(t),
\end{equation}
where $\bbx_i$ is a variable kept at node $i$ and $\bbd_i(t)$ is a local descent direction for node $i$ that depends on iterates $\bbx_{n_i}(t)$. The idea to determine $\bbd_i(t)$ in the decentralized (D)-BFGS method is to let nodes locally approximate the curvature of their own cost functions and those of their neighbors with a local Hessian inverse approximation. We use an update similar to \eqref{eq_bfgs_update} that maintains the secant condition while allowing for decentralized computation. 

To construct such update, define the diagonal normalization matrix $\bbD \in \reals^{np}$ whose $i$th block is $m_i^{-1} \bbI$ and a (small) scalar regularization parameter $\gamma > 0$. Recalling the neighborhood subscript notation, define the modified neighborhood variable and gradient variations, $\tbv_{n_i}(t) \in \reals^{m_i p}$ and $\tbr_{n_i}(t) \in \reals^{m_i p}$, as
\begin{align}
\tbv_{n_i}(t) &:= \bbD_{n_i} \left[ \bbx_{n_i}(t+1) - \bbx_{n_i}(t) \right] \label{eq_dbfgs_vars} \\
\tbr_{n_i}(t) &:= \bbg_{n_i}(t+1) - \bbg_{n_i}(t) - \gamma\tbv_{n_i}(t)\label{eq_dbfgs_grads}.
\end{align}
The neighborhood variations in \eqref{eq_dbfgs_vars} and \eqref{eq_dbfgs_grads} are not simply local components of \eqref{eq_bfgs_vars}. The variable variation in \eqref{eq_dbfgs_vars} differs from the one in \eqref{eq_bfgs_vars} by the presence of the normalizing matrix $\bbD_{n_i}$ and the gradient variation in \eqref{eq_dbfgs_grads} differs by the presence of the term $\gamma\tbv_{n_i}(t)$. Since $\tbv_{n_i}(t)$ and $\tbr_{n_i}(t)$ use only information node $i$ can locally access through neighbors, we can compute and maintain a local Hessian approximation $\bbB^i(t) \in \reals^{m_ip \times m_ip}$, which is updated as the solution of a local regularized version of \eqref{eq_bfgs_update},
\begin{alignat}{2}\label{eq_dbfgs_update}
\bbB^i(t+1) := &\argmin_{\bbZ}\ &&
\text{tr}[ (\bbB^i(t))^{-1} (\bbZ - \gamma \bbI)] - \\\nonumber &&&\qquad\quad
 \text{logdet}[(\bbB^i(t))^{-1} (\bbZ - \gamma \bbI)] - n\\\nonumber
&  \text{   s.t.}\quad &&
\bbZ \tbv_{n_i}(t) = \bbr_{n_i}(t),\quad \bbZ \succeq \bb0.
\end{alignat}
Two properties differentiate \eqref{eq_dbfgs_update} from \eqref{eq_bfgs_update}: (i) The log-determinant forces $\bbB^i(t+1)$ to have eigenvalues greater than $\gamma$. (ii) The secant condition is expressed with respect to the neighborhood gradient variation $\bbr_{n_i}(t)$ and the \emph{modified} neighborhood variable variation $\tbv_{n_i}(t)$. Property (i) is a regularization of \eqref{eq_bfgs_update} first proposed in the context of stochastic quasi-Newton methods \cite{mokhtari2014res}. Property (ii), while not obvious, ensures the secant condition is satisfied as shown in Proposition \ref{prop_global_descent}.

%
\begin{algorithm}[t] 
\setstretch{1.35}
{\small\begin{algorithmic}[1]
  \REQUIRE $\bbB^i(0), \bbx_{i}(0),\bbg_i(0)$, $\bbx_{n_i}(0)$, $\bbg_{n_i}(0)$
  \FOR{$t = 0,1,2, \hdots$}
     \STATE Descent from \eqref{eq_direction_local}:
            $\bbe^i_{n_i}(t) 
                 = - (\bbB^{i}(t)^{-1} +\Gamma\bbD_{n_i}) \bbg_{n_i}(t)$
     \STATE Exchange descent $\bbe^i_j(t) $ with neighbors $j \in n_i$
     \STATE Local descent from \eqref{eq_direction_coord}: 
            $\bbd_{i}(t) :=\sum_{j \in n_i}  \bbe^j_{i}(t).$
     \STATE Local update from \eqref{eq_descent_update_local}:
            $\bbx_{i}(t+1) = \bbx_{i}(t) +\epsilon(t) \bbd_i(t)$ 
     \STATE Exchange $\bbx_{i}(t+1)$ with neighbors. Build $\bbx_{n_i}(t+1)$
     \STATE Compute $\bbg_{i}(t+1) = \nabla_i f(\bbx_{n_i})$ 
            using property  \eqref{eq_gen_problem}
     \STATE Exchange $\bbg_{i}(t+1)$ with neighbors. Build $\bbg_{n_i}(t+1)$
     \STATE Compute $\tbv_{n_i}(t),\tbr_{n_i}(t),\bbB^i(t+1)$ using
            \eqref{eq_dbfgs_vars}, \eqref{eq_dbfgs_grads}, and \eqref{eq_bfgs_dist}
  \ENDFOR
\end{algorithmic}}
\caption{D-BFGS method at node $i$}
\label{alg_dbfgs}
\end{algorithm}

%
For the problem in \eqref{eq_dbfgs_update} to have a solution, we must have $\tbv_{n_i}(t)^T\tbr_{n_i}(t) > 0$ -- see Remark \ref{rmk_inner_product_negative}. When this condition is satisfied, the problem is not only solvable but admits the closed form solution \cite[Proposition 1]{mokhtari2014res}
\begin{align}\label{eq_bfgs_dist}
\bbB^i(t\!+\!1) =  \bbB^i(t)\!+\! \frac{\tbr_{n_i}(t) \tbr^T_{n_i}(t)}{\tbr^T_{n_i}(t) \tbv_{n_i}(t)}   
\!-\! \frac{\bbB^i(t) \tbv_{n_i}(t) \tbv^T_{n_i}(t) \bbB^i(t)}{\tbv^T_{n_i}(t) \bbB^i(t) \tbv_{n_i}(t)} \!+\! \gamma \bbI. 
\end{align}
The differences between \eqref{eq_dbfgs_update} and \eqref{eq_bfgs_update} lead to corresponding differences between \eqref{eq_bfgs_dist} and \eqref{eq_bfgs_cent}. These differences are the addition of the $\gamma \bbI$ term and the use of the variations in \eqref{eq_dbfgs_vars} and \eqref{eq_dbfgs_grads} which are not simple local decompositions of the variations in \eqref{eq_bfgs_vars}. 

\begin{figure}\centering

\def \thisplotscale {0.6}
\def \unit {\thisplotscale cm}

\tikzstyle {block}        = [draw, very thin,
                             rectangle, 
                             minimum height = 1*\unit,
                             minimum width  = 1*\unit]                       

\tikzstyle {blue block}   = [block,
                             fill = blue!20]

\tikzstyle {green  block} = [block,
                             fill = green!20]

\tikzstyle {red block}    = [block,
                             fill = red!20]

\def \deltalabel{ 0.5}

{\fontsize{7}{7}\selectfont\begin{tikzpicture}[x = 1*\unit, y=1*\unit, 
                          shorten >=2pt, shorten <=2pt]

%
%
\path (0,0)        node [red block       ] (Bine) {};
\path (Bine.south) node [red block, below] (Bise) {};
\path (Bise.west)  node [red block, left ] (Bisw) {};
\path (Bisw.north) node [red block, above] (Binw) {};
\path (Bine.north west) ++ (0, \deltalabel) node {$(\bbB^i)^{-1}$}; 
%
%
\path (Binw.west) ++ (-1,0)   node [red   block, left] (rnin) {$\bbg_i$};
\path (Bisw.west) ++ (-1,0)   node [blue  block, left] (rnis) {$\bbg_j$};
\path (rnin.west) ++ (-0.3,0) node [red   block, left] (vnin) {$\bbx_i$};
\path (rnis.west) ++ (-0.3,0) node [blue  block, left] (vnis) {$\bbx_j$};
\path (rnin.north) ++ (0, \deltalabel) node {$\bbg_{n_i}$}; 
\path (vnin.north) ++ (0, \deltalabel) node {$\bbx_{n_i}$}; 
%
%
\path (Bine.east) ++ (1,0) node [red   block, right] (gnin) {$\bbg_i$};
\path (Bise.east) ++ (1,0) node [blue  block, right] (gnis) {$\bbg_j$};
\path (gnin.north) ++ (0, \deltalabel) node {$\bbg_{n_i}$}; 
%
%
\path (gnin.east) ++ (1,0) node [red   block, right] (dnin) {$\bbe_i^i$};
\path (gnis.east) ++ (1,0) node [blue  block, right] (dnis) {$\bbe_j^i$};
\path (dnin.north) ++ (0, \deltalabel) node {$\bbe^i_{n_i}$}; 
%
%
\path (dnin.south east) ++ (1,0) node [red block, right] (dil) {$\bbe_i^i$};
\path (dil.east)        ++ (1,0) node [red block, right] (dir) {$\bbe_i^j$};
\path (rnin.south) -- (Binw.south)  node [midway] {$\text{\eqref{eq_bfgs_dist}} \atop
                                                    \Rightarrow$};
\path (Bine.south) -- (gnin.south)  node [midway] {$\times$};
\path (gnin.south) -- (dnin.south)  node [midway] {$=$};
\path (dil.east)   -- (dir.west)    node [midway, draw, circle, inner sep=1] (plus) {$+$};
\path[draw, -stealth] (plus) -- ++(0,1.5) node [red block, above] {$\bbd_i$};;

%
%
%
\path (0, -3.3)    node [blue block       ] (Bjne) {};
\path (Bjne.south) node [blue block, below] (Bjse) {};
\path (Bjse.west)  node [blue block, left ] (Bjsw) {};
\path (Bjsw.north) node [blue block, above] (Bjnw) {};
\path (Bjse.south west) ++ (0,-\deltalabel) node {$(\bbB^j)^{-1}$}; 
%
%
\path (Bjnw.west) ++ (-1,0)   node [red   block, left] (rnjn) {$\bbg_i$};
\path (Bjsw.west) ++ (-1,0)   node [blue  block, left] (rnjs) {$\bbg_j$};
\path (rnjn.west) ++ (-0.3,0) node [red   block, left] (vnjn) {$\bbx_i$};
\path (rnjs.west) ++ (-0.3,0) node [blue  block, left] (vnjs) {$\bbx_j$};
\path (rnjs.south) ++ (0,-\deltalabel) node {$\bbg_{n_j}$}; 
\path (vnjs.south) ++ (0,-\deltalabel) node {$\bbx_{n_j}$}; 
%
%
\path (Bjne.east) ++ (1,0) node [red   block, right] (gnjn) {$\bbg_i$};
\path (Bjse.east) ++ (1,0) node [blue  block, right] (gnjs) {$\bbg_j$};
\path (gnjs.south) ++ (0,-\deltalabel) node {$\bbg_{n_j}$}; 
%
%
\path (gnjn.east) ++ (1,0) node [red   block, right] (dnjn) {$\bbe_i^j$};
\path (gnjs.east) ++ (1,0) node [blue  block, right] (dnjs) {$\bbe_j^j$};
\path (dnjs.south) ++ (0,-\deltalabel) node {$\bbe^j_{n_j}$}; 
%
%
\path (dnjn.south east) ++ (1,0) node [blue block, right] (djl) {$\bbe_j^j$};
\path (djl.east)        ++ (1,0) node [blue block, right] (djr) {$\bbe_j^i$};
%
%
\path (rnjn.south) -- (Bjnw.south)  node [midway] {$\text{\eqref{eq_bfgs_dist}} \atop
                                                    \Rightarrow$};
\path (Bjne.south) -- (gnjn.south)  node [midway] {$\times$};
\path (gnjn.south) -- (dnjn.south)  node [midway] {$=$};
\path (djl.east)   -- (djr.west)    node [midway, draw, circle, inner sep=1] (plus) {$+$};
\path[draw, -stealth] (plus) -- ++(0,-1.5) node [blue block, below] {$\bbd_j$};;

%
%
\path (dnin.east) edge [-stealth, black!99, bend left]  (dil.north);
\path (dnis.east) edge [-stealth, black!99, bend left]  (djr.north);
\path (dnjs.east) edge [-stealth, black!99, bend right] (djl.south);
\path (dnjn.east) edge [-stealth, black!99, bend right] (dir.south);

%
%
%
\path (rnin.east) edge [-stealth, black!99, bend left] (rnjn.east);
\path (rnjs.west) edge [-stealth, black!99, bend left] (rnis.west);
\path (vnin.east) edge [-stealth, black!99, bend left] (vnjn.east);
\path (vnjs.west) edge [-stealth, black!99, bend left] (vnis.west);

\end{tikzpicture}}
\caption{D-BFGS variable flow. Nodes exchange variable and gradients -- $\bbx_i$ and $\bbg_i$ sent to $j$ and $\bbx_j$ and $\bbg_j$ sent to $i$ -- to build variable and gradient variations $\tbv$ and $\tbr$ that they use to determine local curvature matrices -- $\bbB^i$ and $\bbB^j$. They then use gradients $\bbg_{n_i}$ and $\bbg_{n_j}$ to compute descent directions $\bbe_{n_i}^i$ and  $\bbe_{n_j}^j$. These contain a piece to add locally -- $\bbe_i^i$ stays at node $i$ and $\bbe_j^j$ stays at node -- and a piece to add at neighbors -- $\bbe_j^i$ is sent to node $j$ and $\bbe_i^j$ is sent to node $i$.}
\label{fig_variable_flow_diagram} \end{figure}
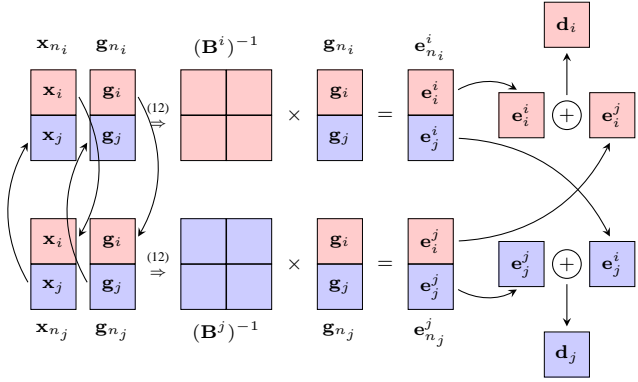

The matrices $\bbB^i$ along with an additional (and also small) regularization parameter $\Gamma >0$ are used by node $i$ to compute the neighborhood descent direction $\bbe^i_{n_i}(t) \in \reals^{m_i p}$ as
\begin{equation}
\bbe^i_{n_i}(t) = - \left( \bbB^i(t)^{-1} + \Gamma \bbD_{n_i} \right) \bbg_{n_i}(t).
\label{eq_direction_local}
\end{equation}
The neighborhood descent direction $\bbe^i_{n_i}(t) \in \reals^{m_ip}$ contains components for variables of node $i$ itself and all neighbors $j \in n_i$ -- see Fig. \ref{fig_variable_flow_diagram}. Likewise, neighboring nodes $j \in n_i$ contain a descent component of the form $\bbe^j_{i}(t)$. The local descent $\bbd_{i}(t)$ is then given by the sum of the components $\bbe^j_{i}(t)$ for all neighbors $j\in n_i$,
\begin{equation}
   \bbd_{i}(t) = \sum_{j \in n_i}  \bbe^j_{i}(t).
\label{eq_direction_coord}
\end{equation}
The descent direction in \eqref{eq_direction_coord} substituted in \eqref{eq_descent_update_local} yields the D-BFGS algorithm outlined in Algorithm \ref{alg_dbfgs}. Each node begins with an initial variable $\bbx_i(0)$, Hessian approximation $\bbB^i(0)$, and gradient $\bbg_i(0)$. Nodes exchange initial variables and gradients to construct initial neighborhood variables $\bbx_{n_i}(0)$ and gradients $\bbg_{n_i}(0)$. For each step $t$, nodes compute their neighborhood descent direction $\bbe_{n_i}(t) $ in Step 2 and exchange the descent elements $\bbe^i_{j}(t)$ with their neighbors in Step 3 to compute the local descent direction $\bbd_i(t)$ in Step 4. They use the local descent direction $\bbd_i(t) $ to update the variable $\bbx_i(t+1)$ and exchange it with their neighbors to form $\bbx_{n_i}(t+1)$ in Steps 5 and 6, respectively. They use these neighbor variables $\bbx_j(t+1)$ to compute an updated local gradient $\bbg_{i}(t+1)$ as in Step 7 and exchange their values in Step 8. In Step 9, nodes compute their neighborhood variable and gradient variations $\tbv_{n_i}(t)$ and $\tbr_{n_i}(t)$ that are required for computing the updated neighborhood Hessian approximation matrix $\bbB^i(t+1)$. 

%

%
An alternative representation of Algorithm \ref{alg_dbfgs} is given in Fig. \ref{fig_variable_flow_diagram} where we emphasize the flow of variables among neighbors. Variable and gradient variations are exchanged -- $\tbv_i(t)$ and $\tbr_i(t)$ are sent to node $j$ and $\tbv_j(t)$ and $\tbr_j(t)$ are sent to node $i$ -- and \eqref{eq_bfgs_dist} is used to compute the curvature estimation matrices -- $\bbB^i(t)$ at node $i$ and $\bbB^j(t)$ at node $j$. The inverses of these matrices are used to premultiply the neighborhood gradients $\bbg_{n_i(t)}$ and $\bbg_{n_j(t)}$, which necessitates an exchange of local gradients -- $\bbg_i(t)$ is sent to node $j$ and $\bbg_j(t)$ to node $i$. This operation results in the computation of the neighborhood descent directions -- $\bbe_{n_i}^i(t)$ and  $\bbe_{n_j}^j(t)$. These descent directions contain a piece to be added locally -- $\bbe_i^i(t)$ stays at node $i$ and $\bbe_j^j(t)$ stays at node -- and a piece to be added at the neighboring node -- $\bbe_j^i(t)$ is sent to node $j$ and $\bbe_i^j(t)$ is sent to node $i$. The local descent direction $\bbd_{i}(t)$ is the addition of the locally computed $\bbe_i^i(t)$ and the remotely computed $\bbe_i^j(t)$ as stated in \eqref{eq_direction_coord}.

%
\subsection{Secant condition in D-BFGS}

To explain the rationale of selecting $\bbB^i(t)$ as in \eqref{eq_dbfgs_update} we show in the following proposition that definitions have been made so that D-BFGS satisfies the secant condition from centralized BFGS.
%
%
\begin{proposition} \label{prop_global_descent}
Consider the D-BFGS method defined by \eqref{eq_descent_update_local}, \eqref{eq_direction_local}, and \eqref{eq_direction_coord} with matrices $\bbB^i(t)$ as given in \eqref{eq_bfgs_dist}. Recall the notational conventions $\bbx(t) = [\bbx_1(t); \ldots; \bbx_n(t)]$ and $\bbg(t)=[\bbg_1(t); \ldots; \bbg_n(t)]$ as well as the definitions of the variable and gradient variations in \eqref{eq_bfgs_vars}. We can rewrite \eqref{eq_descent_update_local}, \eqref{eq_direction_local}, and \eqref{eq_direction_coord} as
\begin{equation}\label{eqn_global_secant_condition}
   \bbx(t+1) = \bbx(t) - \epsilon(t) \big[ \bbH(t) + \Gamma\bbI \big] \bbg(t),
\end{equation}
with a matrix $\bbH(t)$ that satisfies the global secant condition $\bbv(t-1) = \bbH(t) \bbr(t-1)$. 
\end{proposition}

%
\begin{myproof}
Define the matrix $\bbH^i(t) \in \reals^{np \times np}$ to be a block sparse matrix with the sparsity pattern of $n_i$, with a dense sub-matrix $\bbB^{i}(t)^{-1}$, i.e. $\bbH_{n_i}^i(t) = \bbB^{i}(t)^{-1}$. Recall that $\hbx_{n_i}(t) \in \reals^{np}$ stands in for the neighborhood vector $\bbx_{n_i}(t)$ padded with zeros in locations corresponding to non-neighbors. Further recall the matrix $\hbD_{n_i} \in \reals^{np \times np}$ amounts to the matrix $\bbD_{n_i}$ padded with zeros in locations corresponding to non-neighbors. It is apparent then that the global formulation of the neighborhood descent computed by node $i$ from \eqref{eq_direction_local} is $\hbe^i_{n_i}(t) =  - [\bbH^i(t) + \Gamma \hbD_{n_i}] \bbg(t)$. Then, summing over all nodes we have full concatenated descent 
\begin{align}
\bbd(t) = - \sum_{i=1}^n [\bbH^i(t) + \Gamma \hbD_{n_i}] \bbg(t) = - [ \bbH(t) + \Gamma \bbI] \bbg(t),
\end{align}
where we define $\bbH(t) := \sum_{i=1}^n \bbH^i(t)$. To see that $\bbH(t)$ satisfies the secant condition,  Consider then that, by virtue of containing the inverse of a solution to \eqref{eq_dbfgs_update}, $\bbH^i$ satisfies the local secant relation relation $\hbD_{n_i} \bbv(t-1) = \bbH^i(t) \bbr(t-1)$. Again summing over all nodes, we have
\begin{align}
\sum_{i=1}^n  \hbD_{n_i} \bbv(t-1) = \sum_{i=1}^n  \bbH^i(t) \bbr(t-1) ,
\end{align}
which implies the secant condition $\bbv(t-1) = \bbH(t) \bbr(t-1)$.
\end{myproof}

%
The result in Proposition \ref{prop_global_descent} explains the choices in the formulation of the problem that determines the selection of the $\bbB^i(t+1)$ matrices in \eqref{eq_dbfgs_update}. These matrices are chosen so that the overall update in \eqref{eqn_global_secant_condition} satisfies the same secant condition satisfied by (centralized) BFGS.

To clarify the role of the regularization parameters $\gamma$ and $ \Gamma$ we point out that if $\bbB^i(t)$ is positive semidefinite, the constants $\gamma$ and $ \Gamma$ impose the following property on the descent matrix,
\begin{align}
\frac{\Gamma}{ \bar{m}_i}\bbI \preceq \bbB^i(t+1)^{-1} + \Gamma \bbD_{n_i} \preceq \left(\frac{1}{\gamma} + \frac{\Gamma}{\check{m}_i}\right) \bbI  \label{eq_eigen_bounds},
\end{align}
where $\check{m}_i = \min_{j \in n_i} m_j$ and $\bar{m}_i = \max_{j \in n_i} m_j$. In particular, \eqref{eq_eigen_bounds} implies that $\bbB^i(t+1)^{-1}$ is positive semidefinite. Thus, if $\bbB^i(0)^{-1}$ is positive semidefinite, the property in \eqref{eq_eigen_bounds} holds for all times $t$. Inspection of  \eqref{eq_eigen_bounds} shows that the role of $\Gamma$ is to prevent the algorithm from stalling if the eigenvalues of $\bbB^i(t)^{-1}$ become too small. The role of $\gamma$ is to prevent the eigenvalues of $\bbB^i(t)^{-1}$ to become too large. Observe that since it is $\bbB^i(t)^{-1}$ the one that premultiplies $\bbg_{n_i}(t)$, simply adding a regularization factor to \eqref{eq_bfgs_cent} -- which one could mistakenly assume is what we do in \eqref{eq_bfgs_dist} -- may result in a matrix that is very far from satisfying the secant condition. The update in \eqref{eq_bfgs_dist} utilizes the modified gradient and variable variations to pre-compensate for the addition of the $\gamma\bbI$ term so that the secant condition is satisfied {\it after} adding this term. The bounds in \eqref{eq_eigen_bounds} are required for the convergence analyses in Section \ref{sec_convergence}.

%
\begin{remark}\label{rmk_inner_product_negative}\normalfont For the problem in \eqref{eq_dbfgs_update} to have a solution and the update in \eqref{eq_bfgs_dist} to be valid the inner product between the neighborhood variations must be $\tbv_{n_i}(t)^T\tbr_{n_i}(t) > 0$. This condition imposes a restriction in functions that can be handled by D-BFGS. In practical implementations, however, we can check the value of this inner product and proceed to update $\bbB^i(t)$ only when it satisfies $\tbv_{n_i}(t)^T\tbr_{n_i}(t) > 0$.  \end{remark}

%
\section{Asynchronous D-BFGS} \label{sec_async_dbfgs}
Given the amount of coordination between nodes required to implement D-BFGS in Algorithm \ref{alg_dbfgs}, we consider now the D-BFGS algorithm in the asynchronous setting, in which nodes' communications are uncoordinated with those of their neighbors. Our model for asynchronicity follows that used in \cite{bertsekas1989parallel}. Consider that the time indices are partitioned finely enough so that node $i$'s primary computation, namely the computation of descent direction $\bbe^i_{n_i}(t)$, requires multiple consecutive time iterates to complete. For each node $i$, we then define a set $T^i \subseteq \mathbb{Z}^+$ of all time indices in which node $i$ is available to send and receive information, i.e. not busy performing a computation. 

We further define for each node $i$ a function that, given a time index $t$, returns the most recent time node $i$ was available, which we call  $\pi^i(t)$ and expressly define as 
\begin{equation}
\pi^i(t) := \text{max} \{\hat{t} | \hat{t} < t, \hat{t} \in T^i\} \label{eq_time_function}.
\end{equation}
Moreover, we define a function $\pi^i_j(t)$ that, given a time index $t$, returns the most recent time node $j$ sent information that has been received by node $i$ by time $t$, or explicitly, 
\begin{equation}
 \pi_j^i(t) := \pi^j(\pi^i(t)) \label{eq_time_function_neighbor}.
 \end{equation}
In the asynchronous setting, the superscript notation used to denote locally stored information now additionally signifies a node's dated knowledge of a variable,
 \begin{align}
 \bbx_{j}^i(t) := \bbx_j(\pi_j^i(t)), \qquad \bbx_{n_i}^i(t) =[ \bbx_j^i(t)]_{j \in n_i}.
 \label{eq_local_var_async}
 \end{align}
 It is clear then that $\bbx_{j}^i(t) \neq \bbx_j^k(t)$ for any two nodes $i$ and $k$ at any time $t$. We consider as the current global variable state $\bbx(t)$ the concatenation of each node's current knowledge of its own variable, i.e. $\bbx(t) := [\bbx^i_i(t); \hdots; \bbx^n_n(t)]$. We use the same notation for local gradients $\bbg^i_j(t)$ and descent directions $\bbe^i_j(t)$.

 We assume at any time $t \in T^i$ that node $i$ has finished computing a local descent direction it does three things: (i) It reads the variable, gradient, and descent directions from neighboring nodes $j \in n_i$ sent while it was busy. (ii) It updates its local variables and gradient using the descent direction is has just finished computing as well as the descent directions it has received from its neighbors. (iii) Node $i$ can send its locally computed descent direction as well as its updated variable and gradient info. To state in more explicit terms, node $i$ performs the following update to its own block coordinate at all times $t$: 
\begin{align}\label{eq_update_local_async}
    \bbx^i_{i}(t+1) = \bbx^i_{i}(t) + \epsilon(t) \bbd_{i}(t),
\end{align} 
where $\bbd_{i}(t)$ is the decent for the $i$th block $\bbx_{i}(t)$ at time $t$,
\begin{align}
 \bbd_{i}(t) &=  \begin{cases}  \sum_{j \in n_i}  \bbe^j_{i}(t)  & \text{if $t \in T^i$} \\
       \bb0. & \text{otherwise}.
       \end{cases}
       \label{eq_descent_local}
 \end{align}
If $t \in T^i$, node $i$ applies all descent directions available, otherwise it does nothing. Observe that the descent direction in \eqref{eq_descent_local} contains descents calculated with information from time $\pi^i(t)$ and times $\pi^j(t)$ that neighbor $j$ most recently updated its local variable. 
 
To specify the asynchronous version of the decentralized regularized BFGS algorithm, we first reformulate the variable and gradient differences, $\tbv_{n_i}^i(t)$ and $\tbr_{n_i}^i(t)$ for the asynchronous case:
\begin{align}
\tbv_{n_i}^i(t) &= \bbD_{n_i} \left[ \bbx_{n_i}^i(t+1) - \bbx_{n_i}^i(t) \right], \label{eq_dbfgs_vars_async} \\
\tbr_{n_i}^i(t) &= \bbg_{n_i}^i(t+1) - \bbg_{n_i}^i(t) - \gamma \bbv_{n_i}^i(t). \label{eq_dbfgs_grads_async}
\end{align}
We stress that---recalling the superscript notation defined in \eqref{eq_local_var_async}---$\bbx_{n_i}^i(t)$ is the variable state known to $i$ at time $\pi^i(t)$, or the last time node $i$ was available. With this redefined notation, the computation of the local asynchronous BFGS update matrix $\bbB^i(t)$ and the corresponding descent direction $\bbe^i_{n_i}(t)$ follows respectively \eqref{eq_bfgs_dist} and \eqref{eq_direction_local} exactly as in the synchronous setting. 

The complete asynchronous algorithm is outlined in Algorithm \ref{alg_async_dbfgs}. Each node begins with an initial variable $\bbx_i(0)$, Hessian approximation $\bbB^i(0)$, gradient $\bbg_i(0)$, and descent component $\bbe^i_i(0)$. At each time index $t$, they begin by reading the variables of neighbors $\bbe^j_{i}(t), \bbx^i_{j}(t), \bbg^j_{i}(t)$ in Step 2 and construct neighborhood variables. The aggregated descent direction $\bbd_i(t) $ is used to update variables $\bbx_i(t+1)$ and $\bbg_i(t+1)$ in Step 3. Then, with the updated local variable $\bbx_i(t+1)$ and gradient $\bbg_i(t+1)$, node $i$ computes the D-BFGS variables $\tbv^i_{n_i}(t)$, $\tbr^i_{n_i}(t)$, and $\bbB^i(t+1)$ in Step 4. In Step 5, it computes the next descent direction $\bbd^i_{n_i}(t+1)$, and sends its variables to neighbors in Step 6. 


%
\begin{algorithm}[t]
\setstretch{1.35}
{\small\begin{algorithmic}[1]
  \REQUIRE $\bbB^i(0)$, 
           $\bbx_{i}(0)$, 
           $\bbg_{i}(0)$, 
           $\bbe^i_{n_i}(0)$ [cf. \eqref{eq_direction_local}]
  \FOR{$t \in T^i$}
  \STATE Read $\bbe^j_{i}(t), \bbx^i_{j}(t), \bbg^i_{j}(t)$ 
         for $j \in n_i$ from local memory
  \STATE Update $\bbx_{i}(t+1), \bbg_{i}(t+1)$ [cf. \eqref{eq_update_local_async}, \eqref{eq_descent_local}]
  \STATE Compute $\tbv^i_{n_i}(t),\tbr^i_{n_i}(t),\bbB^i(t+1)$ [cf. \eqref{eq_dbfgs_vars_async}, \eqref{eq_dbfgs_grads_async}, \eqref{eq_bfgs_dist}]
  \STATE Compute $\bbe^i_{n_i}(t+1)$ [cf. \eqref{eq_direction_local}]
  \STATE Send $\bbx_{i}(t+1), \bbg_{i}(t+1)$, $\bbe^i_{j}(t+1)$ to neighbors $j \in n_i$  

  \ENDFOR
\end{algorithmic}}
\caption{Asynchronous D-BFGS method at node $i$}
\label{alg_async_dbfgs}
\end{algorithm}

While Algorithm \ref{alg_async_dbfgs} is similar in its basic structure to the synchronous Algorithm \ref{alg_dbfgs}, we highlight a particular difference. In the synchronous algorithm, three rounds of communication were required at each iteration of Algorithm \ref{alg_dbfgs} to properly communicate the dual variable, primal variable, and dual gradient information. In the asynchronous setting, naturally only a single round of communication is possible at each time iteration. As such, all coordination is removed from the algorithm and the order of computation is rearranged slightly in Algorithm \ref{alg_async_dbfgs}. 

 As in the synchronous case, we provide a global formulation of the local descents for the aid in subsequent analysis. While \eqref{eq_descent_local} is an accurate physical description of how the descent is performed by node $i$, the asynchronous setup of \eqref{eq_descent_local} makes it difficult to formulate an equivalent descent direction for the global variable $\bbx(t)$. We alternatively define a virtual formulation for the global descent direction $\bbd(t) \in \reals^{np}$ that is algorithmically equivalent to the one in \eqref{eq_descent_local}, i.e. leads to the same result. Consider the following virtual global update at time $t$, 
  \begin{align}
 \bbx(t+1) = \bbx(t) + \epsilon(t) \bbd(t),
 \label{eq_update_global_async}
 \end{align}
 where the descent direction $\bbd(t) = [ \bbd_{1}(t); \hdots;  \bbd_{n}(t)]$ at time $t$ is
\begin{align}
   \bbd(t) &= \sum_{k:t \in T^k} \hbe_{n_k}^k(t).  \label{eq_descent_global} 
\end{align}
In \eqref{eq_descent_global}, we perform a descent using all directions $\hbe_{n_k}^k(t)$ that finished being computed at time $t$. This is effectively equivalent to assuming that when node $k$ finishes computing a descent direction $\bbe^k_{n_k}(t)$ at time $t \in T^k$, it is instantaneously applied to all neighboring nodes, regardless of whether or not they are busy at time $t$. We assert that, although node $i$ does not physically descend with $\bbe_i^k$ at any time $t \notin T^i$, the virtual update produces the same result as in the physical update for node $i$ at all times $t \in T^i$. This is stated in the following proposition: 
  
\begin{proposition}
The virtual update described by \eqref{eq_update_global_async} and \eqref{eq_descent_global} leads to the same result as the local physical update described by \eqref{eq_update_local_async} and \eqref{eq_descent_local} performed by all nodes.
\label{prop_virtual_update}
\end{proposition}

\begin{myproof}
To show the virtual update is equivalent to the physical update for each node, we first present the coordinate-wise formulation of the virtual update \eqref{eq_descent_global} at time $t$ for node $i$:
   \begin{align}
   \bbd_{i}(t) &=   \sum_{k:t \in T^k} \bbe^k_{i}(t). \label{eq_descent_local_new}
 \end{align}

Consider two nodes $\bbx_i$ and $\bar{\bbx}_i$ who at time $t \in T^i$ are equivalent, i.e. $\bbx_{i}(t) = \bar{\bbx}_{i}(t)$, and will descend asynchronously from their neighbors. Because descent directions calculated by node $i$ are only calculated using information available at times $t \in T^i$, it suffices to show that $\bbx_{i}(t) = \bar{\bbx}_{i}(t)$ at all future $t_+ \in T^i$. 

At time $t$ both nodes compute $\bbe^i_{n_i}(t)$. Node $\bbx_i$ uses \eqref{eq_descent_local} to descent in next iterates while node $\bar{\bbx}_i$ uses \eqref{eq_descent_local_new}. Consider the update performed by the first node at the next available time $t_+$:
\begin{equation}
\bbx_{i}(t_+) = \bbx_{i}(t) + \epsilon \bigg( \sum_{j \in n_i}  \bbe^j_{i}(t) \bigg).
\end{equation}
 Meanwhile, the second node adds descent components as it receives them for all times between $t$ and $t_+$. At time $t_+$, the cumulative update performed by the second node is
 \begin{align}
\bar{\bbx}_{i}(t_+) &= \bar{\bbx}_{i}(t) + \epsilon \sum_{k:t+1 \in T^k} \bbe^k_{i}(t) + \hdots + \epsilon \sum_{k:t_+ \in T^k} \bbe^k_{i}(t) \nonumber \\
&= \bbx_{i}(t) + \epsilon \bigg( \sum_{j \in n_i}  \bbe^j_{i}(t) \bigg) = \bbx_{i}(t_+).
\end{align}
 As this is true for any node $i$, we can also say that the full variable state $\bbx(t_+) = \bar{\bbx}(t_+)$. Furthermore, if it is true that $\bbx(t_+) = \bar{\bbx}(t_+)$, then this will remain the case for all future times $t \in T^i$. 
\end{myproof} 

With Proposition \ref{prop_virtual_update} we show that the global virtual update is equivalent to the physical local update. We continue by establishing the convergence properties of D-BFGS in in both the synchronous and asynchronous settings.

\section{Convergence Analysis}\label{sec_convergence}

We analyze the convergence of D-BFGS method performed on the distributed optimization problem in \eqref{eq_gen_problem} with objective function $f(\bbx)$ with gradient components of the form $\bbg_j(\bbx) = \bbg_j(\bbx_{n_j})$. To begin, we make the following assumption on the eigenvalues of the objective function Hessian,
 \begin{assumption} \label{as_convex}
The objective function $f(\bbx)$ is twice differentiable and the eigenvalues of the objective function Hessian are nonnegative and bounded from above by a positive constant $0 < L < \infty $, 
\begin{equation}
\bb0 \preceq \nabla^2 f(\bbx) \preceq L\bbI.
\label{eq_hessian_bounds}
\end{equation}
\end{assumption} 
 Assumption \ref{as_convex} ensures the objective function $f$ is convex. The upper bound $L$ on the eigenvalues of the Hessian implies that the associated gradient $\bbg(\bbx)$ is Lipschitz continuous with parameter $L$, i.e. $\| \bbg(\bbx) - \bbg(\bbx ')\| \leq  L\|\bbx - \bbx'\|$. In some instances, we can be sure that the $f(\bbx)$ is not just convex, but strongly convex. In these cases, we can show stronger convergence properties of D-BFGS. We therefore introduce the following second assumption.
 \begin{assumption} \label{as_strongly_convex}
The objective function $f(\bbx)$ is twice differentiable and the eigenvalues of the objective function Hessian are nonnegative and bounded from above and below by positive constants $0 < \mu < L < \infty $, i.e. 
\begin{equation}
\mu \bbI \preceq \nabla^2 f(\bbx) \preceq L\bbI.
\label{eq_hessian_bounds_sc}
\end{equation}
\end{assumption} 
In addition to Lipschitz continuity, objective functions that satisfy Assumption \ref{as_strongly_convex} are strongly convex with constant $\mu$. As we will show in Section \ref{sec_consensus} it is possible to derive distributed objective functions for a common class of problems that are both convex and strongly convex. 

We finally make an assumption regarding the inner product of neighborhood variable and gradient variations.
 \begin{assumption} \label{as_inner_product}
For all $i$ and $t$, the inner product between the neighborhood modified variable and gradient vector variations is strictly positive, i.e. $\tbv_{n_i}^T \tbr_{n_i} > 0$.
\end{assumption} 
This assumption is necessary to ensure all local Hessian approximations are well defined in \eqref{eq_bfgs_dist}. While this assumption does not always hold in practice, we use it regardless to simplify analysis. We stress that, in the case the assumption is violated, setting $\bbB^i(t+1) = \bbB^i(t)$ (See Remark \ref{rmk_inner_product_negative}) does not have any bearing on the proceeding analysis.
We proceed to establish the convergence properties of D-BFGS in convex and strongly convex cases.

 \subsection{Synchronous convergence}\label{sec_synch_convergence}

 To discuss the convergence properties of the D-BFGS method in the synchronous setting, we recall that, as established in Proposition \ref{prop_global_descent}, the global descent of D-BFGS can be formulated as $\bbx(t+1) = \bbx(t) - \epsilon(t)[ \bbH(t) + \Gamma \bbI] \bbg(t)$, where $\bbH(t)$ is a matrix built from the local Hessian inverse approximation of each node. The following lemma establishes the positive definiteness of $\bbH^i(t)$ for all $i$ and $t$ with specific bounds on its eigenvalues.
 
\begin{lemma}\label{lemma1}
Consider the D-BFGS method introduced in \eqref{eq_bfgs_dist}-\eqref{eq_direction_coord}. Further, recall both the positive constants $\gamma$ and $\Gamma$ as the regularization parameters of D-BFGS and the definition of the global Hessian inverse approximation $\bbH(t) = \sum_{i=1}^n \bbH^i(t)$. The eigenvalues of the global regularized Hessian inverse approximation $\bbH(t) + \Gamma \bbI$ are uniformly bounded as 
\begin{equation}
\Gamma \bbI \preceq \bbH(t) + \Gamma \bbI \preceq  \Delta \bbI,
\label{eq_prop_eigen_bounds} 
\end{equation} 
where $\Delta := \left( \Gamma + n/\gamma \right)$ and $n$ is the size of network. 
\end{lemma}
 
\begin{myproof}
The lower bound on $\bbH(t) + \Gamma \bbI$ follows immediately from the fact that $\bbH(t)$ is a sum of positive semidefinite matrices and is therefore a positive semidefinite matrix with eigenvalues greater than or equal to 0. The upper bound subsequently follows from the fact that each $\bbH_i(t)$ have eigenvalues upper bounded by $1/\gamma$, as the dense submatrix $\bbB^i(t)^{-1} \preceq 1/\gamma \bbI$. Then, the sum of $n$ such matrices recovers the upper bound in \eqref{eq_prop_eigen_bounds}.
 \end{myproof} 
 
In Lemma \ref{lemma1} we show that there exists lower and upper bounds on the eigenvalues of the Hessian inverse approximation matrix $\bbH(t) + \Gamma \bbI$. From here, it is natural to demonstrate the convergence of the D-BFGS method in the case of either convex or strongly convex functions $f(\bbx)$. In the former case, we show sub-linear convergence of the order of $o(1/t)$ in the proceeding theorem.

\begin{theorem}\label{theorem_convergence}
Consider the D-BFGS method introduced in \eqref{eq_descent_update_local}-\eqref{eq_direction_coord}. If Assumptions \ref{as_convex} and \ref{as_inner_product} hold and the stepsize $\eps(t)$ satisfies $\eps(t) < 2\Gamma/(L \Delta^2)$, then the dual objective function error $f(\bbx(t))-f(\bbx^*)$ converges to zero at least in the order of $o(1/t)$, i.e.,
\begin{equation} \label{eq_convergence_result}
 f(\bbx(t)) - f(\bbx^*)\leq o\left(\frac{1}{t}\right).
\end{equation}
\end{theorem}

\begin{myproof} 
See Appendix A.
 \end{myproof} $\newline$%
 
With Theorem \ref{theorem_convergence} we establish the sub-linear convergence of D-BFGS when the objective function is convex but not strongly convex. By adding a lower bound on the eigenvalues of the Hessian, thus implying strong convexity, we can establish linear convergence as we show in the following theorem.
\begin{theorem}\label{theorem_convergence_sc}
Consider the D-BFGS proposed in  defined in \eqref{eq_descent_update_local}-\eqref{eq_direction_coord}. If Assumptions \ref{as_strongly_convex} and \ref{as_inner_product} hold and stepsize is chosen as $\epsilon(t) < 2\Gamma/(L \Delta^2)$, then the sequence of objective function values $f(\bbx(t))$ converges to the optimal value $f(\bbx^*)$ at least linearly with some constant $0 <c < 1$, i.e. 
\begin{align} \label{eq_convergence_sc_eq}
f(\bbx(t)) - f(\bbx^*) \leq c^t \left( f(\bbx(0)) - f(\bbx^*) \right).
\end{align}
\end{theorem}

\begin{myproof}
See Appendix B.
 \end{myproof} 

With Theorem \ref{theorem_convergence_sc} we establish the linear convergence of D-BFGS in the synchronous setting for a strongly convex objective function. Due to strong convexity, the linear convergence of the sequence $f(\bbx(t)) - f(\bbx^*) \rightarrow 0$ implies the linear convergence of the variable $\| \bbx(t) - \bbx^* \| \rightarrow 0$. We proceed by establishing the convergence properties of asynchronous D-BFGS.

\subsection{Asynchronous convergence}\label{sec_asynch_convergence}
  
  To establish the convergence of decentralized BFGS in the asynchronous setting, it is first necessary to assume a limit to the partial asynchronicity between the nodes. 
  
\begin{assumption}
There exists an asynchronicity limit $0 < B < \infty$ such that, for all $i$, $j$, and $t$, 
\begin{align}
\max\{0, t-B+1\} \leq \pi^i_j(t) \leq t.
\label{eq_async_limit_2}
\end{align}

\label{as_async_limit}
\end{assumption} 
Assumption \ref{as_async_limit} implies a number of things. First, a node available at time $t$ will be available again at least by the time $t+B$. Additionally, any nodes is at most $B$ time iterations out of sync, i.e. node $i$'s knowledge of $\bbx_j$ is at most $B$ descent steps away from the true state of $\bbx_j$. We further assume that a node's communication delay with any other node is bounded by $B$. There are also important implications regarding the convergence of the physical variable update in \eqref{eq_update_local_async} and \eqref{eq_descent_local} with respect to the convergence of the virtual update in \eqref{eq_update_global_async} and \eqref{eq_descent_global}. Specifically, if the the virtual update has converged by time $t^*$, any and all node's local variables will be locally convergent by time $t^* + B$. It is thus sufficient for us to show convergence properties for the virtual update in \eqref{eq_update_global_async}. We proceed to show that the asynchronous D-BFGS algorithm converges with the following theorem.
  
\begin{theorem} \label{theorem_async}
Consider the asynchronous D-BFGS method proposed in \eqref{eq_update_local_async}-\eqref{eq_dbfgs_grads_async} and \eqref{eq_bfgs_dist}-\eqref{eq_direction_local} where $\bbx(0) = \bbx_0$. If Assumptions \ref{as_convex}, \ref{as_inner_product}, and \ref{as_async_limit} hold, then there exists a stepsize $\epsilon(t) > 0$ such that $\lim_{t \rightarrow \infty} \bbg(t) = 0$.
\end{theorem}

\begin{myproof}
See Appendix C.
\end{myproof}  $\newline$
 With the preceding theorem we demonstrate that in the asynchronous setting the the D-BFGS method will indeed converge to the optimal point as time goes to infinity. 
 
We now establish a linear rate of convergence of asynchronous D-BFGS, the rate for synchronous D-BFGS, on a strongly convex function. For the remaining asynchronous analysis we adjust our definition of the asynchronous algorithm slightly to ease the analysis. Given that the discrete time indeces we assign is of our own construction to model real-world time, we can say without loss of generality that only at each time $t$, exactly one node $k$ executes its descent direction $\bbe_{n_k}^k(t)$, i.e. a single term in \eqref{eq_descent_global} rather than a sum. This is equivalent to the time being discretized finely enough so that no two nodes complete the computation of the descent direction at the same time.

To begin, we use an idea used in analysis of incremental gradient algorithms \cite{gurbuzbalaban2015convergence} and first establish a bound on the error between the the asynchronous gradient used by active node $k$, $\bbg^k_{n_k}(t)$, and the $k$th neighborhood component of the true gradient $\bbg_{n_k}(t)$. This is stated formally in the following lemma.
 
\begin{lemma} \label{lemma_error_bound}
Consider the asynchronous D-BFGS algorithm proposed in \eqref{eq_update_local_async}-\eqref{eq_dbfgs_grads_async} and \eqref{eq_bfgs_dist}-\eqref{eq_direction_local}. If Assumptions \ref{as_strongly_convex}, \ref{as_inner_product}, and \ref{as_async_limit} hold, then the norm of the gradient error $\bbdelta_{n_k}(t) := \bbg^k_{n_k}(t) - \bbg_{n_k}(t) $ is upper bounded as 
\begin{equation} \label{eq_lemma_error_bound}
\|\bbdelta_{n_k}(t)\|  \leq 3 \epsilon m_k^2 L^2 \Delta B \max_{t-2B \leq l \leq t-1} \|\bbx(l)-\bbx^* \|.
\end{equation}
\end{lemma}

\begin{myproof}
See Appendix D.
\end{myproof}  $\newline$
 
 With this lemma, we establish that the difference between the synchronous and asynchronous gradient at time $t$ has an upper bound that is proportional to the maximum distance between the optimal variable and the previous $t-2B$ variable states. This is important in establishing a linear convergence rate for asynchronous D-BFGS as we show in the proceeding theorem.

\begin{theorem} \label{theorem_async_linear}
Consider the asynchronous D-BFGS algorithm proposed in \eqref{eq_update_local_async}-\eqref{eq_dbfgs_grads_async} and \eqref{eq_bfgs_dist}-\eqref{eq_direction_local}. If Assumptions \ref{as_strongly_convex}, \ref{as_inner_product}, and \ref{as_async_limit} hold, then with proper choice of stepsize $\epsilon(t) > 0$ such that there exits an $0<c<1$ such that the following holds
\begin{align}
f(\bbx(t)) - f(\bbx^*) \leq c^t (f(\bbx(0)) - \bbx^*)).
\end{align}
\end{theorem}

\begin{myproof}
See Appendix E.
\end{myproof} 
 
 In Theorem \ref{theorem_async_linear} we establish a linear convergence rate for asynchronous D-BFGS, thus demonstrating that introducing asynchronicity between neighboring nodes does not introduce any deterioration to the convergence rate. We proceed to show benefits of D-BFGS numerically by first introducing a common distributed optimization problem called consensus optimization.

%
\section{Consensus Optimization} \label{sec_consensus}

A problem that is often solved distributedly is the minimization of the cost function $\sum_{i=1}^n f_i(\tbx)$ where the variable $\tbx\in\reals^{p}$ is common but the functions $f_i: \reals^{p} \rightarrow \reals$ are locally available at node $i$. This problem can be reformulated into problems that have the structure in \eqref{eq_gen_problem}. To do so, introduce local variables $\bbx_i\in\reals^{p}$ and the aggregate variable $\bbx=[\bbx_1;\dots;\bbx_n]\in\reals^{np}$. The minimization of the sum $\sum_{i=1}^n f_i(\tbx)$ can then be replaced by 
\begin{align}\label{eq_primal_problem}
   \bbx^* \ := \ \argmin_{\bbx\in \reals^{np}} 
                       \  \sum_{i=1}^n f_i(\bbx_i) , \quad
               \ \st \  (\bbI-\bbZ)\bbx=\bb0,
\end{align}
where the matrix $\bbZ\in \reals^{np \times np}$ is chosen so that the feasible variables in \eqref{eq_primal_problem} satisfy $\bbx_i=\bbx_j$ for all $i,j$. A customary choice of a matrix $\bbZ$ with this property is to make it the Kronecker product $\bbZ \coloneqq \bbW \otimes \bbI_p$ of a weight matrix $\bbW\in\reals^{n\times n}$ and the identity matrix $\bbI_p\in\reals^{p \times p}$. The elements of the weight matrix are $w_{ij}> 0$ if $(i,j)\in\ccalE$ and $w_{ij}= 0$ otherwise and the weight matrix $\bbW\in \reals^{n\times n}$ is further assumed to satisfy
\begin{equation}\label{weight_matrix_conditions}
   \bbW                     = \bbW^T,      \quad 
   \bbW\mathbf{1}           = \mathbf{1},  \quad 
   \text{null}\{\bbI-\bbW\} = \text{span}\{\mathbf{1}\}.
\end{equation}
Since $\text{null}(\bbI-\bbW)=\text{span}\{\mathbf{1}\}$, it follows that for any vector $\bbx=[\bbx_1;\dots;\bbx_n]\in \reals^{np}$ the relation $(\bbI-\bbZ)\bbx=\bb0$ holds if and only if $\bbx_1=\dots=\bbx_n$. This means that the feasible variables in \eqref{eq_primal_problem} indeed satisfy $\bbx_i=\bbx_j$ for all $i,j$ and that, consequently, the problem in \eqref{eq_primal_problem} is equivalent to the minimization of $\sum_{i=1}^n f_i(\tbx)$. 

The problem in \eqref{eq_primal_problem} does {\it not} have the structure in \eqref{eq_gen_problem}, but it can be transformed into problems with that structure by the use of penalties in the primal domain, or, alternatively, ascending in the dual domain. We explain this in the following two sections.

%
\subsection{Primal domain penalty methods}\label{sec_primal_consensus}

To transform \eqref{eq_primal_problem} into a formulation with the structure in \eqref{eq_gen_problem} we incorporate the constraint as a penalty term to define the problem
\begin{equation}\label{eq_primal_pen_problem}
   \tbx^* =  \argmin_{\bbx\in \reals^{np}} \sum_{i=1}^n f_i(\bbx_i)
                                           +  \frac{1}{2\alpha} \bbx^T (\bbI - \bbZ) \bbx
          := \argmin_{\bbx\in \reals^{np}} \phi(\bbx),
\end{equation}
where $\alpha$ is a given penalty coefficient. 
The term $(1/2)\bbx^T (\bbI - \bbZ) \bbx$ is a quadratic penalty that pushes $\bbx^*$ to the null space of $(\bbI-\bbZ)^{1/2}$. Since the null spaces of $(\bbI-\bbZ)$ and $(\bbI-\bbZ)^{1/2}$ are identical, this means that $\tbx^*$ is pushed towards the feasible space of \eqref{eq_primal_problem} [cf. \eqref{weight_matrix_conditions}]. The difference between the solutions $\tbx^*$ of \eqref{eq_primal_pen_problem} and $\bbx^*$ of \eqref{eq_primal_problem} is of order $\alpha$; see \cite{YuanQing}.

To compute the gradient of the function $\phi(\bbx)$ that we minimize in \eqref{eq_primal_pen_problem}, begin by observing that the gradient of the penalty term is given by $\nabla[\bbx^T(\bbI-\bbZ)\bbx/2] = (\bbI-\bbZ)\bbx$. Since the matrix $\bbZ$ has a block sparsity pattern that matches the sparsity pattern of the graph, the $i$th component of this gradient can be written as $\nabla_i[\bbx^T(\bbI-\bbZ)\bbx] = \bbx_i - \sum_{j\in n_i} w_{ij} \bbx_j$. As the weights $w_{ij}$ sum up to 1 for any given $i$, we can simplify the latter to $\nabla_i[\bbx^T(\bbI-\bbZ)\bbx] = \sum_{j\in n_i} w_{ij} (\bbx_i - \bbx_j)$. Given that the $i$ component of the gradient of the first sum is simply $\nabla_i \phi(\bbx) = \nabla f_i(\bbx)$. Thus, 
\begin{equation}\label{primal_local_gradient}
   \nabla \phi(\bbx)_i =  \nabla f_{i}(\bbx_i) 
	                    + \frac{1}{\alpha} \sum_{j\in n_i}w_{ij}(\bbx_i - \bbx_j).
\end{equation}
The gradients in \eqref{primal_local_gradient} are locally computable if neighbors exchange variables. The corresponding distributed implementation of gradient descent yields DGD \cite{nedic2009}. In our case, \eqref{primal_local_gradient} is a statement of the distributed computability of the gradient required in \eqref{eq_gen_problem}. We use the explicit form in \eqref{primal_local_gradient} to compute the gradients in Step 7 of Algorithm \ref{alg_dbfgs} or in Step 3 of Algorithm \ref{alg_async_dbfgs} if an asynchronous implementation is preferable. This yields synchronous and asynchronous implementations of a quasi-Newton version of DGD. The local estimation of curvature of this quasi-Newton DGD method results in faster convergence -- see Section VII. 

%
\subsection{Dual ascent methods}\label{sec_dual_consensus}

Introduce the dual variable $\bbnu=[\bbnu_1;\dots;\bbnu_n]\in \reals^{np}$ composed of multipliers $\bbnu_{i}\in \reals^{p}$ that are associated with node $i$ and define the Lagrangian of \eqref{eq_primal_problem} as
\begin{equation}\label{eq_lagrangian}
   \ccalL(\bbx,\bbnu) = \sum_{i=1}^n f_i(\bbx_i) + \bbnu^T(\bbI-\bbZ)\bbx.
\end{equation}
Of importance to our discussion are the primal Lagrangian minimizers that we define as $\bbx(\bbnu) := \argmin_\bbx\ccalL(\bbx,\bbnu)$. Since $\bbZ$ has a block sparsity pattern that matches the sparsity pattern of the graph and the weights $w_{ij}$ sum up to 1 for any given $i$, we can write the second term in \eqref{eq_lagrangian} as $\bbnu^T(\bbI-\bbZ)\bbx = \sum_{i, j\in n_i} w_{ij} \bbx_i^T(\bbnu_i-\bbnu_j)$. Using this fact we conclude that the components $\bbx_i(\bbnu)$ of the Lagrangian minimizer $\bbx(\bbnu)$ are
\begin{equation}\label{eq_lagrangian_minimizers}
   \bbx_i(\nu) = \argmin f_i(\bbx_i) + \sum_{j\in n_i} w_{ij} \bbx_i^T(\bbnu_i-\bbnu_j)
\end{equation}
The Lagrangian minimizers in \eqref{eq_lagrangian_minimizers} can be used to define the dual function $\psi(\bbnu):=\ccalL(\bbx(\bbnu),\bbnu)$ and the corresponding dual problem as finding the argument that maximizes the dual function,
\begin{equation}\label{eq_dual_problem}
   \bbnu^* := \argmax_{\bbnu} \psi(\bbnu) 
            = \argmax_{\bbnu} \ccalL(\bbx(\bbnu),\bbnu) .
\end{equation}
The importance of the optimal dual argument in \eqref{eq_dual_problem} is that the optimal primal argument $\bbx^*$ of \eqref{eq_primal_problem} can be recovered from the Lagrangian minimizer $\bbx(\bbnu^*) := \argmin_\bbx\ccalL(\bbx,\bbnu^*)$ if the primal functions $f_i$ are strongly convex. Another important observation is that gradients of the dual function can be computed by evaluating the constraint slack associated with the Lagrangian minimizers. Specifically, it is not difficult to show that $\nabla \psi(\bbnu) = (\bbI-\bbZ)\bbx(\bbnu)$. Given the block sparsity pattern of $\bbZ$, these gradients can be locally computed as
\begin{equation}\label{eq_dual_derivative}
   \nabla_i \psi(\bbnu)\ =\ \bbx_i - \sum_{j\in n_i} w_{ij} \bbx_j 
                       \ =\ \sum_{j\in n_i} w_{ij} (\bbx_i - \bbx_j) .
\end{equation}
Since the gradients in \eqref{eq_dual_derivative} are functions of neighboring variables only, the distributed computability required in \eqref{eq_primal_problem} holds for the maximization of the dual function in \eqref{eq_dual_problem}. We therefore use  \eqref{eq_dual_derivative} to compute the gradients in Step 7 of Algorithm \ref{alg_dbfgs} or in Step 3 of Algorithm \ref{alg_async_dbfgs}. This yields synchronous and asynchronous implementations of a quasi-Newton version of distributed dual ascent -- see Section VII.

\section{Numerical Results} \label{sec_numerical_results}
We provide numerical results of the performance of D-BFGS on the consensus problem for various objective functions and condition numbers. Simulations are initially performed with the following convex quadratic objective function of variable $\bbx \in \reals^p$.
\begin{align}
f(\bbx) := \sum_{i=1}^n \frac{1}{2} \bbx^T \bbA_i \bbx + \bbb_i^T \bbx
\label{eq_simulation_problem}
\end{align}
where $\bbA_i \in \reals^{p \times p}$ and $\bbb_i \in \reals^p$ define the local objective functions available to node $i$. We control the problems condition number by defining the matrices $\bbA_i = \text{diag}\{ \bba_i \}$. For a chosen condition number $10^{\eta}$, $\bba_i$ is a vector with its $p/2$ elements chosen randomly from the interval $[1, 10^1, \hdots, 10^{\eta/2}]$ and its last $p/2$ elements chosen randomly from the interval $[1, 10^{-1}, \hdots, 10^{-\eta/2}]$, resulting in the sum $\sum_{i=1}^n \bbA_i$ having eigenvalues in the range $[n 10^{-\eta/2}, n 10^{\eta/2}]$. For the vectors $\bbb_i$, the elements are chosen uniformly and randomly from the box $[0,1]^p$.
In our simulations we fix the variable dimension $p=4$ and use a $d$-regular cycle for the graph, in which $d$ is an even number and nodes are connected to their $d/2$ nearest neighbors in either direction. The others parameters such as condition number $10^{\eta}$ and and number of nodes $n$ are varied by simulation. The regularization parameters for BFGS are chosen to be $\gamma = 10^{-2}$ and $\Gamma = 10^{-3}$. In all experiments, we choose a constant stepsize and attempt to pick the largest stepsize for which the algorithms are observed to converge. 

We demonstrate results with solving \eqref{eq_simulation_problem} using both the dual and primal formulations in \eqref{eq_dual_problem} and \eqref{eq_primal_pen_problem}, respectively. The true optimal point $\bbx^*$ can be calculated exactly for the quadratic problem in \eqref{eq_simulation_problem} and we evaluate the average error to be
\begin{align}
\text{error}(t) \coloneqq \frac{1}{n} \sum_{i=1}^n \frac{ \|\bbx_i(t) - \bbx^*\|^2}{ \| \bbx^*\|^2}.
\label{eq_primal_avg_error}
\end{align}

Note that in the dual domain, we find $\bbx(t)$ as the Lagrangian maximizer with respect to the corresponding dual variable $\bbnu(t)$.

\subsection{Synchronous algorithms} 

\begin{figure}[t!]
\vspace{-6mm}
\includegraphics[width=.45\textwidth,height=.25\textheight, keepaspectratio]{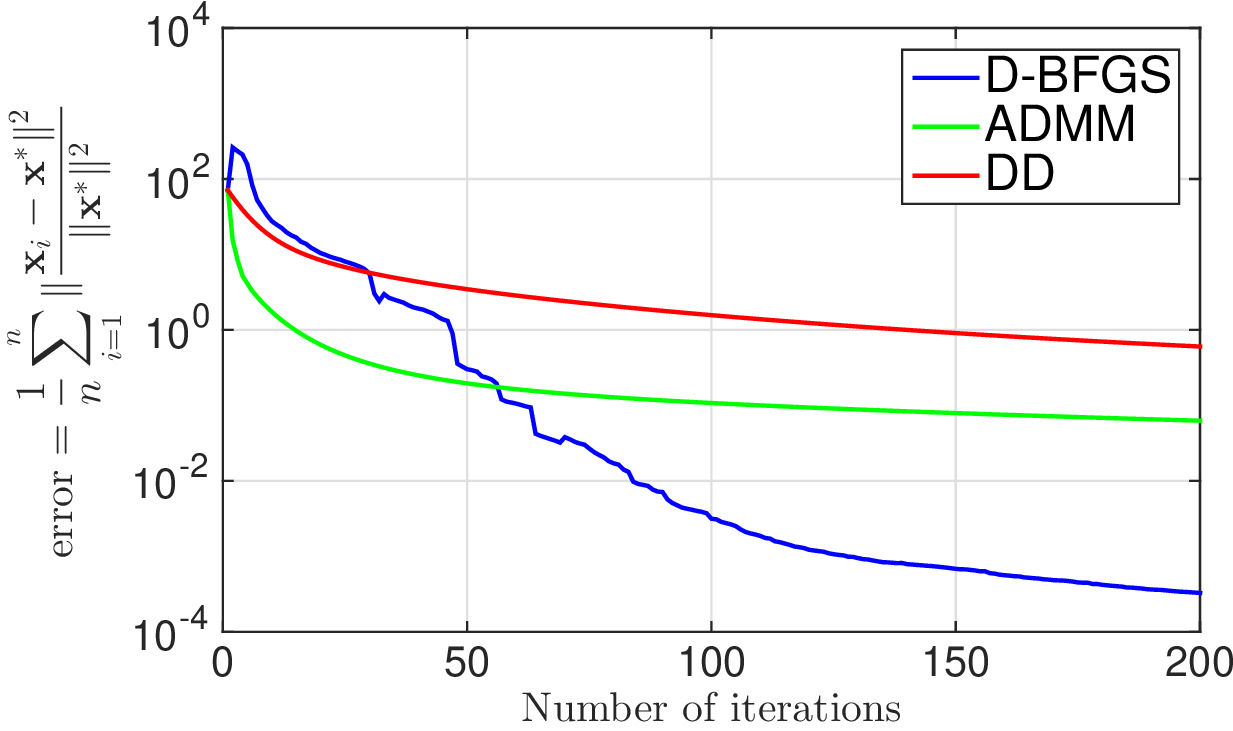}
\vspace{-2mm}
\caption{Convergence of D-BFGS, ADMM, and DD in the dual domain for a quadratic objective function. D-BFGS has the best performance.}
\label{fig_dual_d_4}
\vspace{-2mm}
\end{figure}

\begin{figure}
\includegraphics[width=.48\textwidth,height=.25\textheight]{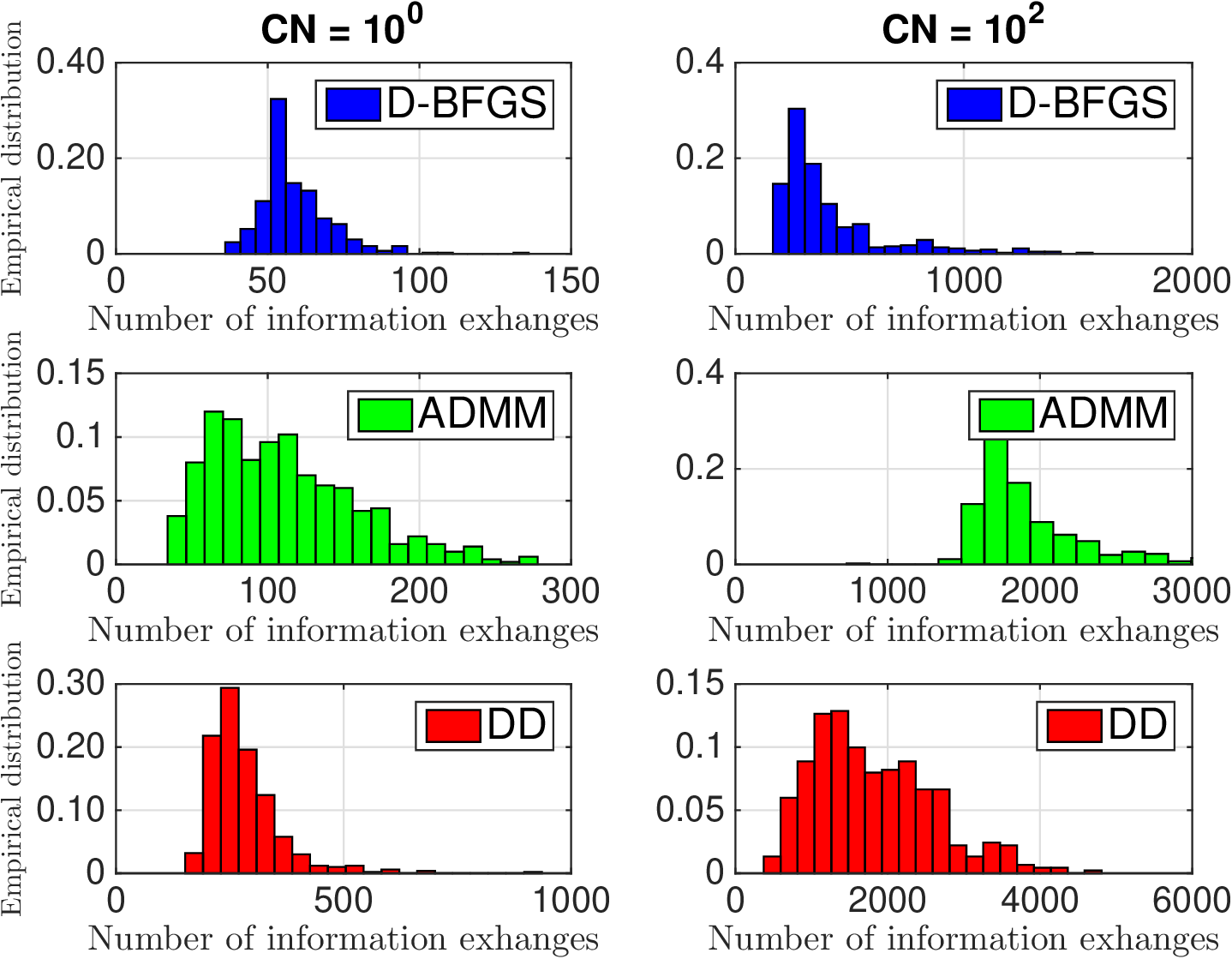}
\caption{Empirical distribution of number of information exchanges needed to reach error of $10^{-2}$ for D-BFGS, ADMM, and DD for quadratic cost function with condition numbers $1$ and $100$. 
For large condition number the gap between the methods is larger.
}
\label{fig_dual_hist}
\end{figure}

We start in by simulating performance in the traditional synchronous setting. We simulate the performance of D-BFGS on the dual problem in \eqref{eq_dual_problem} along with the corresponding first order dual methods, ADMM and DD \cite{BoydEtalADMM11}, using respective stepsizes of $0.01$, $0.002$, and $0.002$ on a network of size $n=50$, connectivity $d=4$, and condition number parameter $\eta=2$. Fig. \ref{fig_dual_d_4} shows the convergence rates of both algorithms in a representative simulation, specifically showing the iteration number vs the average error to the optimal primal variable. Observe that D-BFGS converges substantially faster than both first order methods, achieving an average error of $3\times10^{-4}$ by iteration $200$, while ADMM and DD just reaches average errors of $0.06$ and $0.6$ respectively by iteration 200. We present a more comprehensive view of the difference in convergence times by creating an empirical distribution over a large number of trials. Because D-BFGS requires twice as many communications per iteration as ADMM and DD, Fig. \ref{fig_dual_hist} shows histograms of convergence times of each algorithm in terms of number of local information exchanges. Not only does D-BFGS outperform ADMM and DD in both cases, but the difference in convergence times increases with larger condition number. In particular, there is a factor of 2 between the convergence times of D-BFGS and ADMM for a condition number of $1$ and a factor of 10 for a condition number of $10^2$. The difference in convergence times between D-BFGS and DD is indeed even more significant.

\begin{figure}[t!]
\vspace{-6mm}
\includegraphics[width=.45\textwidth,height=.25\textheight,keepaspectratio]{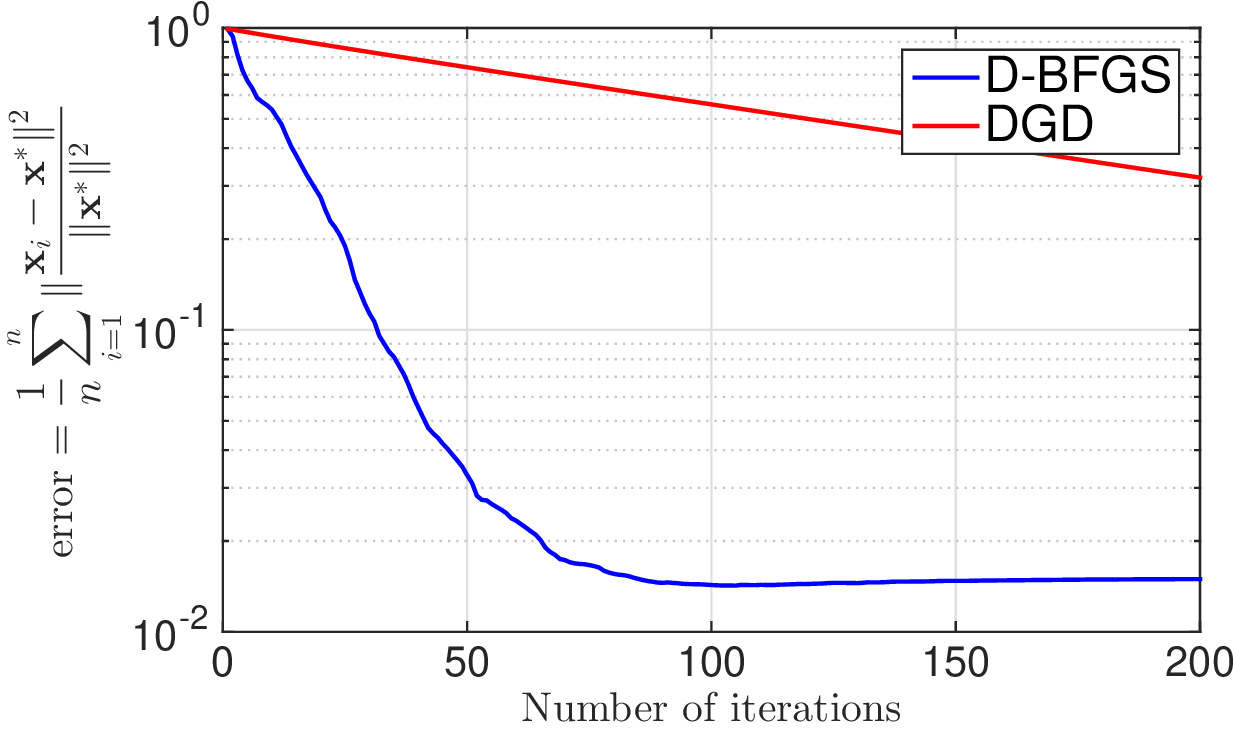}
\vspace{-2mm}
\caption{Convergence of D-BFGS and DGD in the primal domain for a quadratic objective function. D-BFGS converges faster than DGD by almost a factor of 10.}
\label{fig_primal_d_4}
\vspace{-4mm}
\end{figure}

\begin{figure}
\includegraphics[width=.5\textwidth,height=.25\textheight]{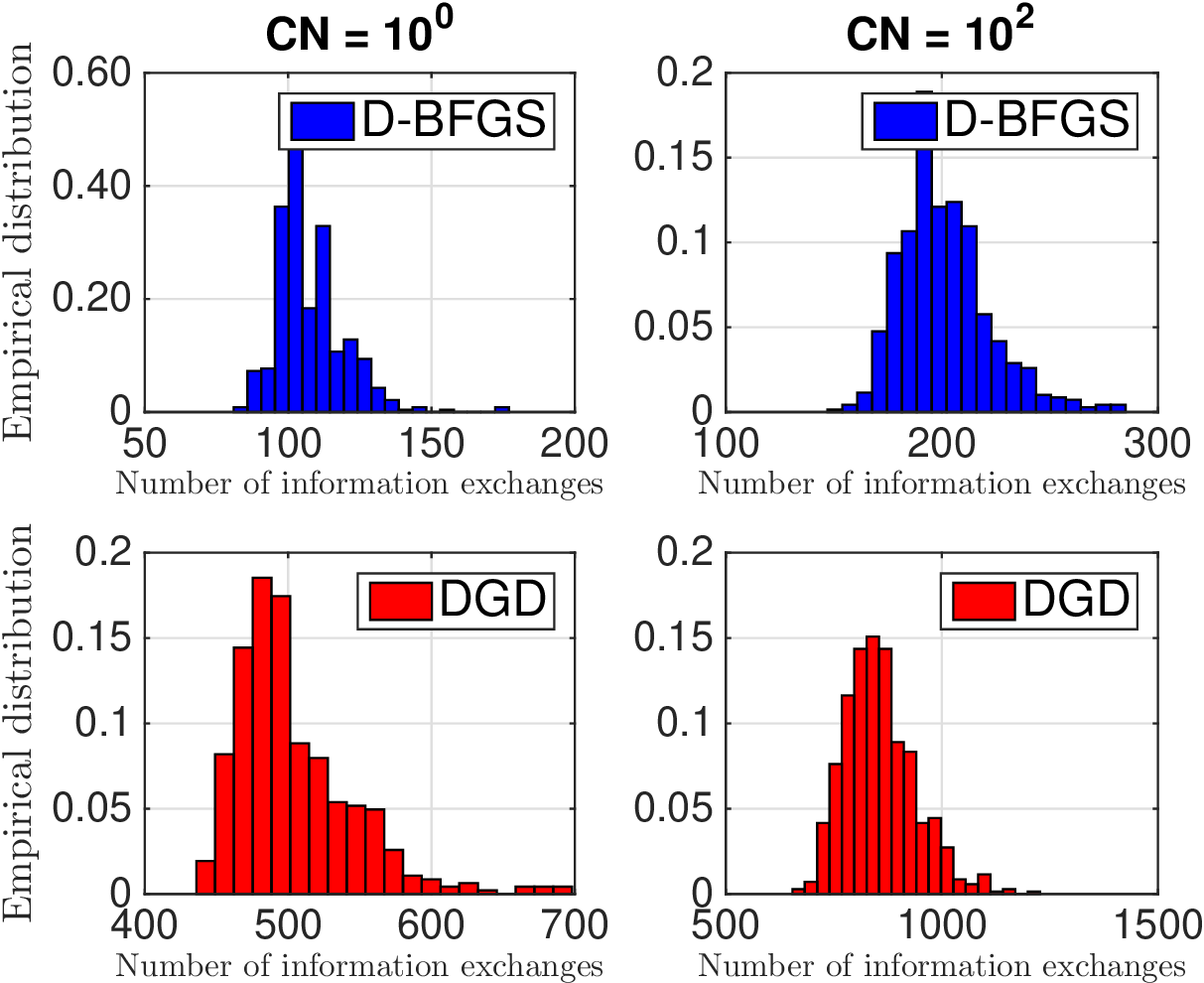}
\caption{Empirical distribution of number of information exchanges needed to reach error of $1.9\times10^{-2}$ for D-BFGS and DGD for quadratic objective function with condition numbers $1$ and $100$. 
The difference in convergence times increases with larger condition numbers.}
\label{fig_primal_hist}
\end{figure}

We numerically evaluate the performance of D-BFGS and DGD on the primal problem in \eqref{eq_primal_pen_problem}, with  parameters set as $n=100$, $d=4$, $p=4$, and $\eta=2$. We additionally set the objective function penalty parameter $\alpha = 10^{-3}$. We choose the row stochastic weight matrix $\bbW$ to be a matrix with diagonal entries $w_{ii} = 1/2 + 1/2(d+1)$ and off diagonal entries $w_{ij} = 1/2(d+1)$ if $j \in n_j$ and $0$ otherwise. The results of a sample simulation, with stepsizes of 0.3 and 1 for D-BFGS and DGD respectively, are shown in Fig. \ref{fig_primal_d_4}. As in the dual domain, D-BFGS converges substantially faster than its gradient descent counterpart, reaching an average error of $0.015$ by iteration 100, while DGD reaches an average error of $0.32$ by iteration 200. A histogram of convergence times with respect to local information exchanges over 1000 realizations is shown in Fig. \ref{fig_primal_hist}. Note in the primal domain D-BFGS requires 3 exchanges per iteration while DGD requires 1. We observe in this case the difference in convergence times between D-BFGS and DGD are around a factor of 5 for both small and large condition numbers.


\subsection{Logistic regression}
To evaluate the performance of D-BFGS on problem of more practical interest, we additionally look at the logistic regression problem. In logistic regression, we seek to learn a linear classifier $\bbx$ that can predict the label of a data point $v_j \in \{-1,1\}$ given a feature vector $\bbu_j \in \reals^p$. To do so, we evaluate for a set of training samples the likelihood of a label given a feature vector as $P(v=1 | \bbu) = 1/(1+\exp(-\bbu^T\bbx))$ and find $\bbx$ that maximizes the log likelihood over all samples. In the distributed setting, it is often assumed that the training set is large and distributed amongst $n$ servers, with server $i$ receiving $q_i$ samples. It is then the case that each server $i$  has access to a different objective function given the training samples $\{\bbu_{il}\}_{l=1}^{q_i}$ and $\{v_{il}\}_{l=1}^{q_i}$. The aggregate objective function can be defined as 
\begin{align}
f(\bbx) := \frac{\lambda}{2}\| \bbx\|^2 + \sum_{i=1}^{n} \sum_{l=1}^{q_i} \log[ 1 + \exp(-v_{il}\bbu_{il}^T\bbx)],
\label{eq_logistic_problem}
\end{align}
where the first term is a regularization term used to reduce overfitting and is parametrized by $\lambda \geq 0$. 

\begin{figure}[t]
\vspace{-7mm}
\includegraphics[width=.45\textwidth,height=.28\textheight, keepaspectratio]{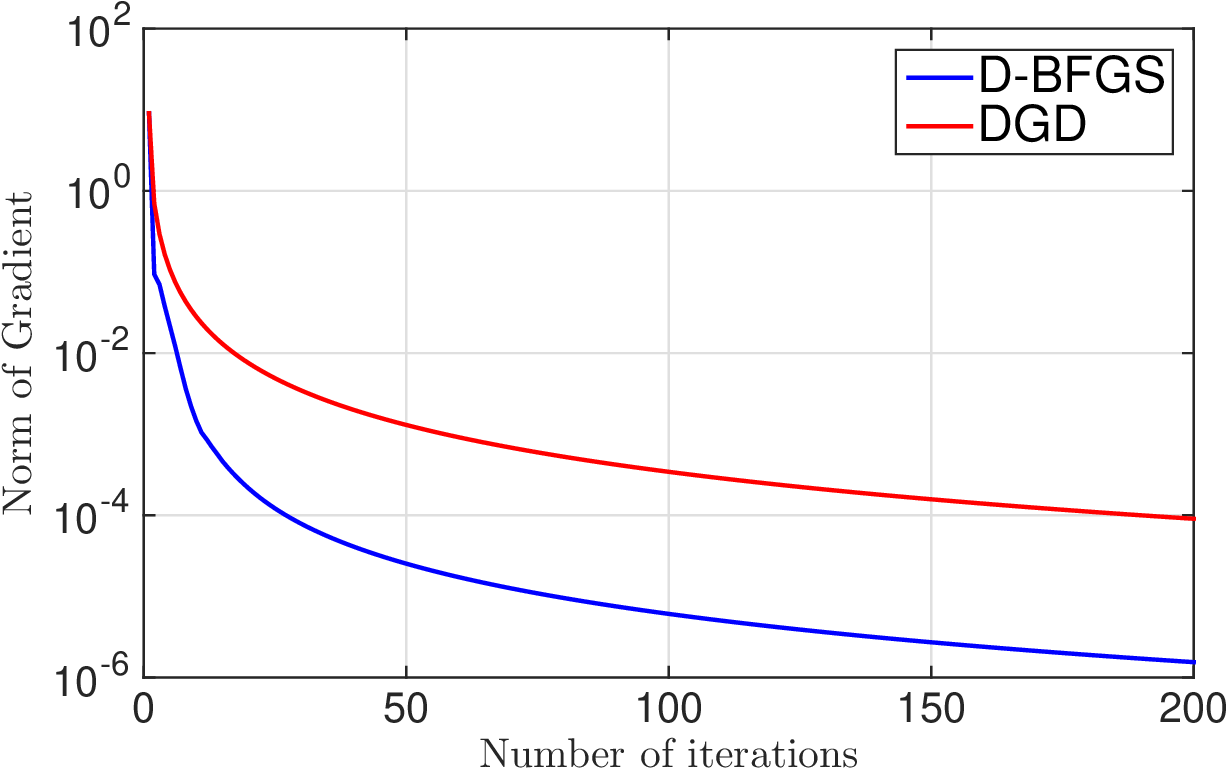}
\caption{Convergence of D-BFGS and DGD in the primal domain for a logistic regression problem.}
\label{fig_log_reg_plot}
\end{figure}

For our simulations we generate a dataset of feature vectors $\bbu_{il}$ with label $v_{il}=1$ from a normal distribution with mean $\mu$ and standard deviation $\sigma_{+}$, and with label $v_{il}=-1$ from a normal distribution with mean $-\mu$ and standard deviation $\sigma_{-}$. Each node $i$ receives $q_i=100$ samples and the regularization parameters is fixed to be $\lambda= 10^{-4}$. The feature vector parameters are set as $\mu=3$ and $\sigma_{+}=\sigma_{-}=1$ to make the data linearly separable.

The other parameters we set the same as in earlier simulations, i.e. $n=100$ nodes connected in $d=4$-regular cycle with $p=4$. The D-BFGS regularization parameters are chosen as $\Gamma = \gamma = 10^{-1}$ with stepsizes of $0.3$ and $1$ for D-BFGS and DGD respectively. The resulting convergence paths are shown in Fig. \ref{fig_log_reg_plot}, in this case shown with respect to the norm of the gradient. D-BFGS reaches a gradient magnitude of $1.3 \times 10^{-6}$ before iteration $200$ with DGD reaching a gradient magnitude of $9.1 \times 10^{-5}$.

\subsection{Asynchronous algorithms}

We compare the performance of D-BFGS and DD in the asynchronous setting on dual formulation of the quadratic problem in \eqref{eq_simulation_problem}. The model we use for the asynchronicity is modeled after a random delay phenomenon in physical communication systems that creates asynchronous local clocks between servers. Each server $i$'s local clock begins at $t^i_0 = t_0$ and selects subsequent times as $t^i_k = t^i_{k-1} + N^i_{k}$, where $N^i_k$ is drawn from a normal distribution with mean $\mu$ and standard deviation $\sigma$. The standard deviation effectively controls the level of drift or asynchronicity between the nodes, i.e. a larger standard deviation will lead to nodes deviating further from the local clocks of their neighbors.

\begin{figure}[t]
\vspace{-7mm}
\includegraphics[width=.45\textwidth,height=.25\textheight,keepaspectratio]{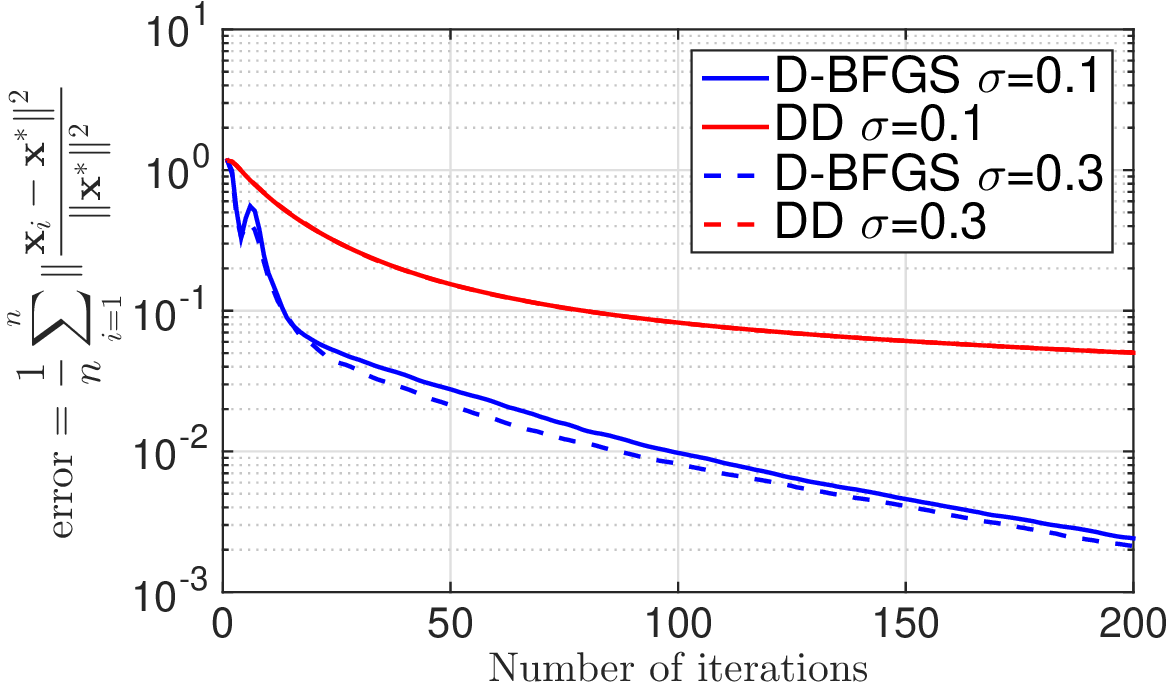}
\caption{Convergence of asynchronous D-BFGS and DD in the dual domain for a quadratic objective function. The level of asynchronicity varies between $\sigma=0.1$ and $\sigma=0.3$.}
\label{fig_a_dual_d_4}
\end{figure}

In our initial experiment, we set $n=50$ nodes with dimension $p=4$, condition number parameter $\eta = 1$, and network connectivity $d=4$. The D-BFGS regularization parameters are set to be $\gamma = \Gamma = 10^{-1}$ and use the same stepsizes used in the synchronous setting ($0.01$ and $0.002$ respectively for D-BFGS and DD). For our asynchronicity parameters we set $\mu=1$ and select two standard deviations $\sigma=0.1$ and $\sigma=0.3$. The resulting convergence paths are shown in Fig. \ref{fig_a_dual_d_4}. For this figure the number of iterations in the x-axis refers to the average error after all nodes have reached that number of iterations locally. Observe that for $\sigma=0.01$ D-BFGS outperforms DD, reaching by iteration 200 and average error of $2.4 \times 10^{-3}$ and DD reaching only an average error of $5.1 \times 10^{-2}$. Further observe that, for the case of $\sigma=0.3$, the larger drift does not result in a substantially different convergence time for either method, suggesting that the performance of the asynchronous algorithms is not very sensitive to changes in the level of asynchronicity between nodes.  


\section{Conclusions} \label{sec_conclusion}
We considered the problem of general decentralized optimization, in which nodes sought to minimize a cost function while only being aware of a local strictly convex component. The problem was solved through the introduction of D-BFGS as a decentralized quasi-Newton method. In D-BFGS, each node approximates the curvature of its local cost function and its neighboring nodes to correct its descent direction. Analytical results were established in both synchronous and asynchronous versions of the algorithm. We also showed numerical results on two types of consensus optimization problems in both the dual and primal domains, in which significant improvement was observed over alternatives.


\section*{Appendix A: Proof of Theorem \ref{theorem_convergence}} \label{sec_theorem_convergence}
The steps of the proof follows closely those of Proposition 1.3.3 in \cite{bertsekas1999nonlinear} for gradient descent for (not strongly) convex functions. Given the Lipschitz continuity of $f$ in \eqref{eq_hessian_bounds}, we have the following inequality for $f(\bbx(t+1))$ for constant stepsize $\epsilon(t) = \epsilon$:
\begin{align}\label{proof0}
f(\bbx(t\!+\!1)) \leq&  f(\bbx(t))\! -\! \epsilon\bbg(t)^T (\bbH(t)+\Gamma\bbI) \bbg(t) \! \\ &\qquad +\! \frac{L\epsilon^2}{2} \|(\bbH(t)+\Gamma\bbI) \bbg(t) \|^2. \nonumber
\end{align}
%
%
Using the lower and upper bounds on the eigenvalues of $\bbH(t) + \Gamma \bbI$ for the second and third term, respectively, we can write
\begin{align}
f(\bbx(t+1)) &\leq f(\bbx(t)) - \epsilon \| \bbg(t) \|^2 \left[ \Gamma - \frac{\epsilon L  \Delta^2}{2} \right ]. \label{proof1}
\end{align}
Assume that $\epsilon$ satisfies $\Gamma - \epsilon L  \Delta^2/2 > 0$. Denote by $\mathcal{N}^*$ the set of global minima and distance between $\bbx(t)$ and the set of minima
\begin{align} \label{eq_distance}
d(\bbx(t),\mathcal{N}^*) := \text{min}_{\bbx^* \in \mathcal{N}^*} \| \bbx(t) - \bbx^* \|
\end{align}
The convexity of $h(\bbx)$ implies that for any $\bbx^* \in \ccalN^*$ we have
$
f(\bbx(t)) \leq f(\bbx^*) + \bbg(t)^T (\bbx(t)- \bbx^*).
$
This inequality in conjunction with the Cauchy-Swartz inequality implies that
\begin{align}
f(\bbx(t)) &\leq  f(\bbx^*) + \|\bbg(t)\| \|(\bbx(t)- \bbx^*\|.
\end{align}
As this holds for all $\bbx^* \in \ccalN^*$, we can obtain
\begin{align}
f(\bbx(t)) - f(\bbx^*) &\leq \| \bbg(t)\| d(\bbx(t),\mathcal{N}^*). \label{proof2}
\end{align}
For notational convenience, define $e(t) := f(\bbx(t)) - f(\bbx^*)$ and assume without loss of generality that $d(\bbx(t),\mathcal{N}^*) \neq 0$. Now we combine the results of \eqref{proof1} and \eqref{proof2} and rearrange terms to get
\begin{align}
e(t+1) &\leq e(t)\left[ 1 - \epsilon \left( \Gamma - \frac{\epsilon L \Delta^2}{2} \right) \frac{e(k)}{d(\bbx(t),\mathcal{N}^*)^2} \right]. \label{proof3}
\end{align}
The inequality in \eqref{proof3} indeed implies that $e(t) \leq o(1/t)$. The details of this derivation are provided in the aforementioned proposition in \cite{bertsekas1999nonlinear}, which we remove for space considerations.

\section*{Appendix B: Proof of Theorem \ref{theorem_convergence_sc}} \label{sec_theorem_convergence_sc}

Consider the D-BFGS update with constant stepsize $\bbx(t+1) = \bbx(t) - \epsilon(\bbH(t) + \Gamma\bbI) \bbg(t)$ along with Taylor's expansion of the function $f$ in \eqref{proof0}.
%
%
%
%
We use the upper and lower bounds on the eigenvalues of $\bbH(t) + \Gamma \bbI$ from Lemma \ref{lemma1} and subtract $f^* := f(\bbx^*)$ from both sides to upper bound \eqref{proof0} as
%
\begin{align}
f(\bbx(t\!+\!1)) \!-\! f^* \!\leq\! f(\bbx(t)) \!-\! f^*
\!-\!\epsilon\left[\Gamma\!-\! \frac{L \Delta^2 \epsilon}{2}\right]\! \|\bbg(t) \|^2.
\end{align}
As a result from strong convexity, we have $\|\bbg(t) \|^2 \geq 2\mu(f(\bbx(t))-f^*)$. If we choose $\epsilon < 2\Gamma/(L \Delta^2)$, we subsequently have after rearranging terms
\begin{align} \label{prooflin}
&f(\bbx(t+1)) - f^*\leq [f(\bbx(t)) - f^*] \left[1 - 2\mu\epsilon\left(\Gamma - \frac{L \Delta^2 \epsilon}{2}\right)\right].
\end{align}
%
%
%
We obtain linear convergence if $0 < 1 - \left(2\mu\Gamma\epsilon -\mu {L \Delta^2 \epsilon^2}\right) < 1$, which holds for out previous choice of  $0< \epsilon < 2\Gamma/(L  \Delta^2)$. Expanding \eqref{prooflin} $t$ times we achieve the result in \eqref{eq_convergence_sc_eq} with $c = 1 - \left(2\mu\Gamma\epsilon -\mu {L \Delta^2 \epsilon^2}\right) $.

%

\section*{Appendix C: Proof of Theorem \ref{theorem_async}} \label{sec_theorem_async}
The steps of our analysis follow closely that of asynchronous gradient descent in \cite[Proposition 5.1]{bertsekas1989parallel}. Consider the global virtual descent formulation in \eqref{eq_update_global_async} and \eqref{eq_descent_global}. As established in Proposition \ref{prop_virtual_update}, this is equivalent to the asynchronous formulation. The Hessian eigenvalue bounds in \eqref{eq_hessian_bounds_sc} allow us to write
\begin{align}
f(\bbx(t+1))  \leq& f(\bbx(t))  + \epsilon \bbd(t)^T \bbg(t) + \frac{L}{2} \epsilon^2 \|\bbd(t)\|^2 \label{eq_descent_lemma}.
\end{align}
We look first at bounding the second term in the summand. Recall the hat notated $\hbg^k(t) \in \reals^{np}$, which signifies the local vector $\bbg^k_{n_k}(t)$ padded with zeros. Further define $\bbH_{n_k}^k(t) := \bbB^k(t)^{-1} + \Gamma\bbI$. We can then substitute for $\bbd(t)$ and rearrange as follows,
\begin{align}
\bbd(t)^T \bbg(t) =& -\bbg(t)^T \sum_{k:t \in T^k} \hbH_{n_k}^k(t) \hbg_{n_k}^k(t) \\
=& - \sum_{k:t \in T^k} \bbg^T_{n_k}(t) \bbH_{n_k}^k(t) \bbg_{n_k}^k(t).
\end{align}
We proceed by adding and subtracting the asynchronous gradient $\bbg_{n_k}^k(t)$ and rearranging terms to obtain
\begin{align}
\bbd(t)^T \bbg(t) =& \sum_{k:t \in T^k} \Big[ -\bbg_{n_k}^{kT}(t)  \bbH_{n_k}^k(t) \bbg_{n_k}^k(t) \label{eq_middle_0}  \\
& \quad + (\bbg_{n_k}^k(t) - \bbg_{n_k}(t))^T \bbH_{n_k}^k(t) \bbg_{n_k}^k(t) \Big]. \nonumber
\end{align}
We bound the first and second terms in summand using the lower and upper eigenvalue bounds in \eqref{eq_eigen_bounds} and the Cauchy-Schwartz inequality, respectively.
\begin{align}\label{eq_middle_1} 
&\bbd(t)^T \bbg(t) \\
&\leq \!\sum_{k:t \in T^k}\! \Big[\! -\Gamma \|\bbg_{n_k}^k(t)\|^2 + \Delta\|\bbg_{n_k}^k(t) - \bbg_{n_k}(t)\| \|\bbg_{n_k}^k(t)\| \Big]  \nonumber.
\end{align}
Next, we bound $\|\bbg_{n_k}^k(t) - \bbg_{n_k}(t)\|$, which represents the difference between the actual gradient and the asynchronous gradient seen by node $k$. Split the norm into its components as
\begin{align}
\|\bbg_{n_k}^k(t) - \bbg_{n_k}(t)\| &\leq \sum_{j \in n_k} \|\bbg_{j}^k(t) - \bbg_{j}(t)\|. \label{eq_middle_2}
\end{align}
 We bound $\|\bbg_{j}^k(t) - \bbg_{j}(t)\|$ by noting that they represent components of the gradient of the global variable at times $t$ and $\pi_j^k(t)$,
\begin{align}
\|\bbg_{j}^k(t) - \bbg_{j}(t)\| &= \| \nabla f ( \bbx( \pi^k_j(t)))_j - \nabla f ( \bbx( t))_j\|.
\end{align}
This can be bounded using the gradients Lipschitz continuity as
\begin{align}
\|\bbg_{j}^k(t) - \bbg_{j}(t)\|&\leq L \| \bbx( \pi^k_j(t))) - \bbx(t)\| \label{eq_middle_3}.
\end{align}
Recall the asynchronicity bound $B$ that limits that amount of time between $\pi^k_j(t)$ and $t$. We proceed in bounding the term on difference between the global and asynchronous variable as 
\vspace{-1mm}
\begin{align}
\| \bbx( \pi^k_j(t))) - \bbx(t)\| &\leq \epsilon \| \sum_{\tau=t-\pi^k_j(t)}^{t-1} \bbd(\tau)\| \leq \epsilon \| \sum_{\tau=t-B}^{t-1} \bbd(\tau)\| \nonumber.
\end{align}
From the triangle inequality we obtain then
\vspace{-1mm}
\begin{align}
\| \bbx( \pi^k_j(t))) - \bbx(t)\| &\leq \epsilon \sum_{\tau = t-B}^{t-1} \| \bbd(\tau) \| \label{eq_variable_delay_bound}.
\end{align}
We can substitute the results in \eqref{eq_variable_delay_bound}, \eqref{eq_middle_3}, and \eqref{eq_middle_2} back into \eqref{eq_middle_1} and rearrange terms to get 
\begin{align}\label{eq_middle_4}
&\bbd(t)^T \bbg(t) \\ 
&\leq \sum_{k:t \in T^k} \Big[ -\Gamma \|\bbg_{n_k}^k(t)\|^2 + \epsilon m_k L \Delta\|\bbg_{n_k}^k(t)\|  \sum_{\tau = t-B}^{t-1} \| \bbd(\tau) \| \Big] .\nonumber
\end{align}
Finally, we can substitute the result in \eqref{eq_middle_4} back into \eqref{eq_descent_lemma} to obtain
\begin{align} \label{eq_descent_lemma_2}
f(\bbx(t+1))  &\leq f(\bbx(t))  +  \frac{L}{2} \epsilon^2 \|\bbd(t)\|^2  
 - \! \epsilon \Gamma \bigg[\sum_{k:t \in T^k} \! \|\bbg_{n_k}^k(t)\|\bigg]^2  \nonumber \\  
&\quad + \epsilon^2   \bar{m} L \Delta  \sum_{k:t \in T^k}  \|\bbg_{n_k}^k(t)\| \sum_{\tau=t-B}^{t-1} \| \bbd(\tau)\|,
\end{align}
where we introduce the term $\bar{m} := \text{max}_k \{ m_k \}$ for notational convenience. Note that the third term on the right hand side was further bounded using the triangle inequality. We simplify notation by introducing the variable
$K(t) \coloneqq \sum_{k:t \in T^k} \|\bbg_{n_k}^k(t)\|.$
We proceed by bounding $\| \bbd(t)\|$ using the triangle inequality as
\begin{align}
\| \bbd(t)\| &= \| \sum_{k:t \in T^k} \hbH_{n_k}^k(t) \hbg_{n_k}^k(t) \| \leq  \sum_{k:t \in T^k} \| \bbH_{n_k}^k(t) \bbg_{n_k}^k(t) \| \nonumber.
\end{align}
Using the upper bound on the eigenvalues of $\bbH^k_{n_k}$ be then obtain
\begin{align}
\| \bbd(t)\|  &\leq \Delta \sum_{k:t \in T^k} \|\bbg_{n_k}^k(t)\| = \Delta K(t). \label{eq_descent_bound} 
\end{align}
Replace $\| \bbd(t)\|$ in \eqref{eq_descent_lemma_2} by the upper bound in \eqref{eq_descent_bound} to obtain
\begin{align}
f(\bbx(t+1))  &\leq f(\bbx(t)) - \left( \Gamma \epsilon - L \Delta^2 \epsilon^2/2 \right) K(t)^2   \nonumber \\  &\qquad + L \Delta^2 \epsilon^2 \bar{m}  K(t) \sum_{\tau=t-B}^{t-1} K(\tau)  \label{eq_descent_lemma_3}.
\end{align}
The last term in \eqref{eq_descent_lemma_3} can be bounded further using the inequality $|a| |b| \leq a^2 + b^2$ and rearranging terms to obtain
\begin{align}
f(\bbx(t+1))  &\leq f(\bbx(t)) - \left[ \Gamma \epsilon - \frac{L \Delta^2 \epsilon^2}{2} - B L \Delta^2 \epsilon^2 \bar{m} \right] K(t)^2   \nonumber \\  & \qquad + L \Delta^2 \epsilon^2  \bar{m} \sum_{\tau=t-B}^{t-1} K(\tau)^2 
 \label{eq_descent_lemma_4}.
\end{align}
By adding the inequality in \eqref{eq_descent_lemma_4} or all $\tau$ between 0 and $t$, notice that the last term includes a summation of $K(\tau)^2$ from $0$ to $t$, with each term repeated $B$ times. We then obtain
\begin{align}
f(\bbx(t+1))  \leq f(\bbx_0) \label{eq_descent_lemma_5}  \!-\! \left[ \Gamma \epsilon - \frac{L \Delta^2 \epsilon^2}{2} - 2 B L \Delta^2 \epsilon^2 \bar{m} \right] \!\sum_{\tau=0}^{t}  K(\tau)^2.  
\end{align}
Supposing we choose $0 < \epsilon  < \Gamma/(L \Delta^2/2 + 2BM \bar{m} \Delta^2)$, then the second term in \eqref{eq_descent_lemma_5} is positive. We subtract optimal value $f^*$ from both sides and, noting that $f(\bbx(t+1))  \geq f^*$, rearrange terms to obtain
%
%
%
\vspace{-1mm}
\begin{align}
\sum_{\tau=0}^{t} K(\tau)^2 \leq 
\frac{f(\bbx_0)  - f^*}{\Gamma \epsilon - \frac{L \Delta^2 \epsilon^2}{2} - 2 B L \Delta^2\epsilon^2 \bar{m} } .  \label{eq_descent_lemma_7}
\end{align}
Following the assumption that $f(\bbx_0)  - f^*$ is bounded and positive, we conclude that the limit of the summand in the left hand side must go to zero,
\begin{align}
\lim_{\tau \to \infty} K(\tau)  = \lim_{\tau \to \infty} \sum_{k:\tau \in T^k} \| \bbg^k(\tau)) \| = 0.
\label{eq_descent_lemma_limit}
\end{align}
We substitute \eqref{eq_descent_lemma_limit} into \eqref{eq_descent_bound} to obtain $\lim_{\tau \to \infty} \| \bbd(\tau) \| = 0$ and by extension with \eqref{eq_variable_delay_bound} that $\lim_{\tau \to \infty} \| \hbx^k(\tau) - \bbx(\tau)\| = 0$ for all $k$. The Lipschitz continuity condition that follows from \eqref{eq_hessian_bounds} yields 
\begin{align}
\lim_{\tau \to \infty}  \| \hbg^k(\tau) - \bbg(\tau)) \| = 0.
\label{eq_gradient_difference_limit}
\end{align}
\vspace{-1mm}
Finally, we conclude from \eqref{eq_descent_lemma_limit} and the partial asynchronicity assumption that $\lim_{\tau \to \infty} \| \hbg^k(\tau) \| = 0$ for all $k$ and, with \eqref{eq_gradient_difference_limit}, we have 
%
$\lim_{\tau \to \infty}  \| \bbg(\tau) \| = 0.
\label{eq_gradient_limit}$
%
Thus, the global virtual variable $\bbx(t)$ and, by partial asynchronicity, all local $\bbx_i(t)$ are convergent.

\section*{Appendix D: Proof of Lemma \ref{lemma_error_bound}} \label{sec_lemma_error_bound}
The steps of this proof are adapted from Section 3.2 in \cite{gurbuzbalaban2015convergence}. We begin to find an upper bound on the norm of $\bbdelta_{n_i}(t)$ by considering the definition along with the bounds from \eqref{eq_middle_2}-\eqref{eq_variable_delay_bound}.
\begin{align}
\|\bbdelta_{n_k}(t)\| &= \| \bbg^k_{n_k}(t) - \bbg_{n_k}(t) \| \leq \epsilon m_k L \sum_{\tau = t-B}^{t-1} \| \bbd(\tau) \|   \nonumber  
\end{align}
We can subsequently bound $\|\bbdelta_{n_k}(t)\|$ using the bound from \eqref{eq_descent_bound}, 
\begin{align}
\|\bbdelta_{n_k}(t)\| &\leq \epsilon m_k L \Delta \sum_{\tau = t-B}^{t-1} \| \bbg^k_{n_k}(\tau) \|.
\end{align}
Using the definition of $\bbdelta(t)$ and the triangle inequality on the final factor and substitute the bound from \eqref{eq_middle_3}, we then have
\begin{align}
\|\bbdelta_{n_k}(t)\| &\leq \epsilon m_k L \Delta \sum_{\tau = t-B}^{t-1} \left[ \| \bbg_{n_k}(\tau) \| + \| \bbe_{n_k}(\tau)\| \right] \nonumber.
\end{align}
We bound the second summand with an alternative bound of $\|\bbdelta_{n_k}(t)\|$. To do so, we use the bound from \eqref{eq_middle_2}-\eqref{eq_middle_3} to obtain
\begin{align}
\|\bbdelta_{n_k}(t)\| &\leq \epsilon m_k L \Delta \sum_{\tau = t-B}^{t-1} [ \| \bbg_{n_k}(\tau) \| \nonumber \\
& \qquad + L \sum_{j \in n_k} \| \bbx(\pi^k_j(\tau)) - \bbx(\tau) \| ].
\end{align}
We proceed in bounding the final term by adding and subtracting $\bbx^*$ and then using the triangle inequality to obtain
\begin{align}
\|\bbdelta_{n_k}(t)\| &\leq \epsilon m_k L \Delta \sum_{\tau = t-B}^{t-1} [ \| \bbg_{n_k}(\tau) \|  \\
& \qquad + L \sum_{j \in n_k} \left(\| \bbx(\pi^k_j(\tau)) - \bbx^*\| + \|\bbx^* - \bbx(\tau) \| \right) ]. \nonumber 
\end{align}
Take the maximum over all time iterations between $\tau-B$ and $\tau$ to bound both first and second term in the final sum to obtain
\begin{align}
\|\bbdelta_{n_k}(t)\| &\leq \epsilon m_k L \Delta \sum_{\tau = t-B}^{t-1} [ \| \bbg_{n_k}(\tau) \| \label{bb2}  \\
& \qquad + 2 m_k L \max_{\tau-B \leq l \leq \tau} \| \bbx(l) - \bbx^*\| ]. \nonumber  
\end{align}
To combine terms in the sum in \eqref{bb2}, consider that we can bound the first summand on the right hand side using Lipschitz continuity and then similarly take the maximum over $\tau-B\leq l \leq\tau$ to obtain
\begin{align}
\|\bbdelta_{n_k}(t)\|&\leq \epsilon m_k L \Delta \sum_{\tau = t-B}^{t-1} 3 m_k L \max_{\tau-B \leq l \leq \tau} \| \bbx(l) - \bbx^*\|.
\end{align}
We obtain our final result in \eqref{eq_lemma_error_bound} by increasing the range of the maximum to include all $t-2B \leq l \leq t-1$ and summing $B$ times
\begin{align}
\| \bbdelta_{n_k}(t) \| &\leq 3 \epsilon m_k^2 M^2 \Delta B \max_{t-2B \leq l \leq t-1} \| \bbx(l) - \bbx^*\|. 
\end{align}
%

\section*{Appendix E: Proof of Theorem \ref{theorem_async_linear}} \label{sec_theorem_async_linear}
Consider the following that results from Lipschitz continuity,
\begin{align}
f(\bbx(t+1)) &\leq f(\bbx(t) - \epsilon \bbg(t)^T \bbd(t) + \frac{L \epsilon^2}{2} \| \bbd(t) \|^2 \nonumber.
\end{align}
We substitute the asynchronous $\bbd(t) = -\hbH_{n_k}^k(t) \hbg_{n_k}^k(t)$--where we again notate by $\bbH_{n_k}^k(t) := \bbB^k(t)^{-1} + \Gamma\bbI$ and $k$ is the active node at time $t$--and add and subtract the true gradient $\bbg(t)$ from the second two terms. After applying the upper eigenvalue bound of $\bbH_{n_k}^k(t)$ on the final term and rearranging terms, we obtain
\begin{align}
f(\bbx(t+1)) &\leq f(\bbx(t)) - \epsilon \Gamma \| \bbg(t) \|^2 + \epsilon \bbg(t)^T \hbH_{n_k}^k(t) (\bbg(t) - \hbg_{n_k}^k(t)) \nonumber \\
& \quad + \epsilon^2 L \Delta^2 \| \bbg(t)\|^2 + \epsilon^2 L \Delta^2 \| \bbg_{n_k}(t) - \bbg^k_{n_k}(t) \|^2  \nonumber. 
\end{align}
We can then apply the Cauchy-Schwartz inequality and the upper eigenvalue bound of $\bbH^k(t)$ to the third term in the previous expression. After rearranging terms we have
\begin{align}
f(\bbx(t+1)) &\leq f(\bbx(t)) - \epsilon \left( \Gamma - \epsilon L \Delta^2 \right) \| \bbg(t) \|^2  \label{bb3} \\ 
& \qquad +  \epsilon \Delta \|\bbg(t)\| \| \bbdelta(t)\| + \epsilon^2 L \Delta^2 \| \bbdelta(t) \|^2. \nonumber 
\end{align}
We substitute the gradient error bound from \eqref{eq_lemma_error_bound} into \eqref{bb3},
\begin{align}
f(\bbx(t+1)) &\leq f(\bbx(t)) - \epsilon \left( \Gamma - \epsilon L \Delta^2 \right) \| \bbg(t) \|^2  \\
& \qquad +  C \epsilon^2 \Delta \|\bbg(t)\| \max_{t-2B \leq l \leq t-1} \| \bbx(l) - \bbx^*\| \nonumber \\
& \qquad + C^2 \epsilon^4 L \Delta^2 (\max_{t-2B \leq l \leq t-1} \| \bbx(l) - \bbx^*\|)^2,  \nonumber 
\end{align}
where we define the constant $C := 3 m_k^2 M^2 \Delta B$ for notational convenience. We use the bound given by strong convexity to bound $\| \bbx(l) - \bbx^*\| \leq \| \bbg(l) \|/\mu$ and combine terms to obtain
\begin{align}
f(\bbx(t+1)) &\leq f(\bbx(t)) - \epsilon \left( \Gamma - \epsilon L \Delta^2 \right) \| \bbg(t) \|^2  \\
& \qquad + \frac{C \epsilon^2 \Delta}{\mu}\left(1 + C \epsilon^2 L \Delta/\mu \right)  \max_{t-2B \leq l \leq t-1} \| \bbg(l)\|^2.  \nonumber 
\end{align}
Subtract $f^* := f(\bbx^*)$ from both sides of the inequality. In addition, we can bound the second and third terms respectively by the common lower and upper bounds on the gradient norm, i.e. $\| \bbg(t) \|^2 \geq 2\mu( f(\bbx(t)) - f^*)$ and $\| \bbg(t) \|^2 \leq 2L( f(\bbx(t)) - f^*)$. After substitution of these bounds, we obtain
\begin{align}\label{bb5}
&f(\bbx(t+1)) - f^* \\
&\leq f(\bbx(t)) - f^* - 2\epsilon \mu \left( \Gamma - \epsilon L \Delta^2 \right) ( f(\bbx(t)) - f^*) \nonumber  \\
& \quad + \frac{2 L C \epsilon^2 \Delta}{\mu}\left(1 + C \epsilon^2 L \Delta/\mu \right) \max_{t-2B \leq l \leq t-1}( f(\bbx(l)) - f^*).  \nonumber 
\end{align}
To establish linear convergence, we repeat from \cite[Lemma 3]{feyzmahdavian2014delayed}.
\begin{lemma}\label{lemma_lin_convergence_async}
Consider the a nonnegative sequence ${V_t}$ and constants $p,q > 0$ satisfying
\begin{equation}
V(t+1) \leq pV(t) + q \max_{t - d(t) \leq l \leq t} V(l).
\end{equation}
If $p + q <1$ and $0 \leq d(t) \leq d_{max}$ for some constant $d_max >0$, then the sequence converges at a linear rate, i.e. 
\begin{equation}
V(t) \leq (p+q)^{\frac{t}{1+d_{max}}} V(0).
\end{equation}
\end{lemma}

We conclude the proof by restating \eqref{bb5} as follows:
\begin{align}
f(\bbx(t+1)) - f^* &\leq p (f(\bbx(t)) - f^*) \label{bb6}\\ 
& \quad + q \max_{t-2B \leq l \leq t-1}( f(\bbx(l)) - f^*),  \nonumber 
\end{align}
where $p = 1 - 2\epsilon \mu(\Gamma - \epsilon L \Delta^2)$ and $q = 2LC\epsilon^2 \Delta (1 + C\epsilon^2L\Delta/\mu)/\mu$. Choosing $\epsilon$ small enough such that $p+q<1$ holds, we then have linear convergence as a result of Lemma \ref{lemma_lin_convergence_async}.


\urlstyle{same}
\bibliographystyle{IEEEtran}
\bibliography{bmc_article.bib}

\end{document}